\PassOptionsToPackage{unicode}{hyperref}
\PassOptionsToPackage{hyphens}{url}
\PassOptionsToPackage{dvipsnames,svgnames,x11names}{xcolor}
\documentclass[
  12pt]{article}
\usepackage{makecell}
\usepackage{amsmath,amssymb}
\usepackage{iftex}
\usepackage{algorithm}
\usepackage{float}
\usepackage{algpseudocode}
\usepackage{multirow}
\ifPDFTeX
  \usepackage[T1]{fontenc}
  \usepackage[utf8]{inputenc}
  \usepackage{textcomp} 
\else 
  \usepackage{unicode-math}
  \defaultfontfeatures{Scale=MatchLowercase}
  \defaultfontfeatures[\rmfamily]{Ligatures=TeX,Scale=1}
\fi
\usepackage{lmodern}
\ifPDFTeX\else  
\fi
\IfFileExists{upquote.sty}{\usepackage{upquote}}{}
\IfFileExists{microtype.sty}{
  \usepackage[]{microtype}
  \UseMicrotypeSet[protrusion]{basicmath} 
}{}
\makeatletter
\@ifundefined{KOMAClassName}{
  \IfFileExists{parskip.sty}{%
    \usepackage{parskip}
  }{
    \setlength{\parindent}{0pt}
    \setlength{\parskip}{6pt plus 2pt minus 1pt}}
}{
  \KOMAoptions{parskip=half}}
\makeatother
\usepackage{xcolor}
\makeatletter
\ifx\paragraph\undefined\else
  \let\oldparagraph\paragraph
  \renewcommand{\paragraph}{
    \@ifstar
      \xxxParagraphStar
      \xxxParagraphNoStar
  }
  \newcommand{\xxxParagraphStar}[1]{\oldparagraph*{#1}\mbox{}}
  \newcommand{\xxxParagraphNoStar}[1]{\oldparagraph{#1}\mbox{}}
\fi
\ifx\subparagraph\undefined\else
  \let\oldsubparagraph\subparagraph
  \renewcommand{\subparagraph}{
    \@ifstar
      \xxxSubParagraphStar
      \xxxSubParagraphNoStar
  }
  \newcommand{\xxxSubParagraphStar}[1]{\oldsubparagraph*{#1}\mbox{}}
  \newcommand{\xxxSubParagraphNoStar}[1]{\oldsubparagraph{#1}\mbox{}}
\fi
\makeatother

\usepackage{longtable,booktabs,array}
\usepackage{calc} 
\usepackage{etoolbox}
\makeatletter
\patchcmd\longtable{\par}{\if@noskipsec\mbox{}\fi\par}{}{}
\makeatother
\IfFileExists{footnotehyper.sty}{\usepackage{footnotehyper}}{\usepackage{footnote}}
\makesavenoteenv{longtable}
\usepackage{graphicx}
\makeatletter
\def\maxwidth{\ifdim\Gin@nat@width>\linewidth\linewidth\else\Gin@nat@width\fi}
\def\maxheight{\ifdim\Gin@nat@height>\textheight\textheight\else\Gin@nat@height\fi}
\makeatother
\setkeys{Gin}{width=\maxwidth,height=\maxheight,keepaspectratio}
\makeatletter
\def\fps@figure{htbp}
\makeatother

\makeatletter
\@ifpackageloaded{caption}{}{\usepackage{caption}}
\AtBeginDocument{%
\ifdefined\contentsname
  \renewcommand*\contentsname{Table of contents}
\else
  \newcommand\contentsname{Table of contents}
\fi
\ifdefined\listfigurename
  \renewcommand*\listfigurename{List of Figures}
\else
  \newcommand\listfigurename{List of Figures}
\fi
\ifdefined\listtablename
  \renewcommand*\listtablename{List of Tables}
\else
  \newcommand\listtablename{List of Tables}
\fi
\ifdefined\figurename
  \renewcommand*\figurename{Figure}
\else
  \newcommand\figurename{Figure}
\fi
\ifdefined\tablename
  \renewcommand*\tablename{Table}
\else
  \newcommand\tablename{Table}
\fi
}
\@ifpackageloaded{float}{}{\usepackage{float}}
\floatstyle{ruled}
\@ifundefined{c@chapter}{\newfloat{codelisting}{h}{lop}}{\newfloat{codelisting}{h}{lop}[chapter]}
\floatname{codelisting}{Listing}

\makeatother
\makeatletter
\@ifpackageloaded{caption}{}{\usepackage{caption}}
\@ifpackageloaded{subcaption}{}{\usepackage{subcaption}}
\makeatother

\ifLuaTeX
  \usepackage{selnolig}  
\fi
\usepackage[]{natbib}
\usepackage{bookmark}

\IfFileExists{xurl.sty}{\usepackage{xurl}}{} 
\hypersetup{
  pdftitle={piecewise Polynomial Probability Density Estimation with Quantile-Based Binning},
  pdfauthor={Author 1; Author 2},
  pdfkeywords={3 to 6 keywords, that do not appear in the title},
  colorlinks=true,
  linkcolor={blue},
  filecolor={Maroon},
  citecolor={Blue},
  urlcolor={Blue},
  pdfcreator={LaTeX via pandoc}}

\newcommand{\anon}{1}

\begin{document}

\def\spacingset#1{\renewcommand{\baselinestretch}%
{#1}\small\normalsize} \spacingset{1}


\if1\anon
{
  \title{\bf Moment-based Piecewise Polynomial Probability Density Estimation with Quantile-Based Binning}
  \author{Meltem Turan$^a$ and
    Joakim Munkhammar$^a$\thanks{
    Corresponding author \\ \textit{Email address: joakim.munkhammar@angstrom.uu.se} (Joakim Munkhammar)}\hspace{.5cm}\\
    $^a$Department of Civil and Industrial Engineering, Uppsala University, \\Uppsala, Sweden}
  \maketitle} \fi

\if0\anon
{
  \bigskip
  \bigskip
  \bigskip
  \begin{center}
    {\LARGE\bf Moment-based Piecewise Polynomial Probability Density Estimation with Quantile-Based Binning}
\end{center}
  \medskip
} \fi

\bigskip
\begin{abstract}
Accurate reconstruction of probability density functions (PDFs) from data is essential in engineering applications. Classical global moment-based polynomial approximations often suffer from oscillations, instability in the tails, and sensitivity to the choice of support. This work proposes a quantile-based piecewise polynomial density reconstruction approach that combines equal-probability binning with local moment-matched polynomials within each bin. Two variants are considered: piecewise monomial and piecewise Lagrange polynomials with Chebyshev nodes. The numbers of bins and polynomial degrees are selected by a proposed grid search approach guided by the Kolmogorov--Smirnov (K--S) test statistic under non-negativity constraints. Across several benchmark distributions, the proposed methods reduce K--S errors by about $80$--$96\%$ relative to standard monomial and Lagrange polynomial approaches, and by about $83$--$97\%$ compared with spline density estimation. For real-world household electricity consumption and solar irradiance data, the piecewise approaches achieve K--S test statistic performance comparable to kernel density estimation while offering improved control over tail behavior and oscillations. Overall, the results demonstrate that quantile-based localization substantially enhances the robustness and fidelity of moment-based polynomial PDF reconstruction.
\end{abstract}

\noindent%
{\it Keywords:} Nonparametric estimation; Method of Moments; piecewise polynomial approximation; Quantile-based binning; Density reconstruction.
\vfill

\newpage
\spacingset{1.8} 

\section{Introduction}\label{sec-intro}

The probability density function (PDF) is a fundamental object in statistics and probability theory, describing how probability is distributed over the support of a random variable. In applications, the underlying PDF is rarely known in closed form and must be estimated from observed data \cite{silverman2018density}. Obtaining accurate and robust PDF estimates is important not only for purely statistical tasks (e.g., inference, goodness-of-fit testing) but also for engineering applications such as forecasting, reliability estimation, and risk assessment \cite{jain2022load,li2022integrated}.

In energy engineering, PDF estimation is especially relevant for modeling stochastic processes such as the temporal variability of solar irradiance \cite{munkhammar2018markov}, household electricity consumption \cite{munkhammar2021very}, and wind power generation \cite{han2019kernel,gavaskar2018fast}. In these contexts, distributions often exhibit skewness, heavy tails, or multimodality, and can deviate substantially from standard parametric families. Estimators must therefore be flexible enough to capture such complex behavior while remaining numerically stable and interpretable.

Broadly, PDF estimation methods can be categorized into parametric and nonparametric methods. Parametric methods represent the data using a prescribed family of distributions, e.g., Normal, Log-Normal, Weibull, or Gamma \cite{shi2021wind}. Parameters are typically estimated via the Method of Moments or maximum likelihood. Moment-based methods, in particular, approximate a distribution by matching low-order moments (mean, variance, skewness, kurtosis, etc.) to those of the observed data \cite{alibrandi2018kernel}. Although this approach can be efficient and interpretable when the chosen family is appropriate, it suffers from model misspecification: when the true distribution does not belong to the assumed family, moment estimates may become biased and the resulting PDF inaccurate \cite{valiant2017estimating}.

Nonparametric methods, by contrast, seek to reconstruct the distribution directly from the data without assuming a specific parametric form \cite{gu2021forecasting}. A widely used nonparametric method is Kernel Density Estimation (KDE), which approximates the density by placing a smooth kernel function at each data point \citep{pavlides2022non}. The crucial hyperparameter in KDE is the bandwidth, which controls the amount of smoothing. A single global bandwidth is typically used across the entire support, which may be suboptimal: dense regions require fine resolution, whereas sparse regions (e.g., tails) benefit from stronger smoothing. This limitation can lead to oversmoothing of sharp features or undersmoothing of noisy regions \citep{reyes2017bandwidth}. Another line of nonparametric methods involves histogram-based approaches and their refinements \citep{Freedman1981}. Quantized Histogram-Based Smoothing (QHBS) is designed to overcome the discontinuities of standard histograms by smoothing bin counts to obtain a continuous and differentiable density estimate \cite{hsu2021moon}. QHBS can capture multimodal or discontinuous distributions, but its performance remains sensitive to the choice of the number of bins and binning strategy \cite{dhal2021histogram}. Piece-wise polynomial-based probability density estimation has mainly been developed in the context of histogram smoothing and orthogonal polynomial-based methods such as B-spline functions \cite{Chan2013}. Early work interprets histograms as piece-wise constant polynomials, with optimal bin-width rules mitigating—but not eliminating—bias–variance trade-offs and discontinuities at bin edges \cite{Abdous2001,reyes2017bandwidth}. Spline Density Estimation (SDE) constructs the PDF as a non-negative mixture of B-spline basis functions, providing a flexible piecewise-polynomial representation with local support and smoothness control \cite{papp2014shape}. One defines a set of B-splines over the data range, typically with clamped knots at the endpoints to ensure coverage of the support. The spline coefficients are then obtained via a non-negative least-squares problem, guaranteeing non-negativity of the estimated PDF. However, the quality of the spline approximation depends heavily on the number and placement of knots: too few basis functions may oversmooth the density, while too many may overfit noise. Fixed spline bases may also struggle to adapt to regions with varying data density or to represent sharp peaks \citep{kirkby2023spline}. Recent studies construct higher-order piece-wise polynomials on finite supports, enforcing non-negativity and smoothness via quadratic programming and moment or interpolation constraints \cite{Yu2023}. Among them, moment-based polynomial approaches approximate the probability density on a fixed interval by polynomials expressed in monomial, Chebyshev, Hermite, or Lagrange polynomial bases, which are orthogonal polynomials, Coefficients are determined by matching a finite number of moments \citep{munkhammar2017polynomial,turan2024polynomial}. These methods can be effective on suitably chosen intervals with an appropriate number of moments, but they often degrade in the presence of small variance, heavy tails, strong skewness, or poorly tuned intervals. Moreover, moment-based polynomial approximations are prone to oscillations (Runge-type behavior), misplacement of modes, and negativity issues that require ad hoc adjustments of the degree or support.

This study presents and evaluates a modified moment-based polynomial approximation approach that introduces local adaptivity via quantile-based binning. Instead of applying a single global polynomial on a fixed interval, the approach partitions the data into bins using empirical quantiles and fits a separate polynomial in each bin via the Method of Moments. This approach can be viewed as a quantile-based piecewise polynomial reconstruction of the PDF. The key idea is to leverage the binning mechanism of QHBS, but rather than smoothing bin counts, the density within each bin is approximated locally using polynomials based on the moments. The bin boundaries are chosen so that each bin contains approximately the same number of observations, thereby addressing data sparsity in tail regions and providing a more balanced representation across the support.

The proposed approach does not aim to produce a perfectly smooth global approximation; instead, it focuses on reducing oscillations and improving local fidelity, especially in tails and around peaks of skewed or multimodal distributions. Two polynomial bases, which are piecewise monomial and piecewise Lagrange polynomials with Chebyshev interpolation points in each bin, are considered for the local approximations. For both variants, the coefficients are obtained by matching local empirical moments in each bin. The number of bins $N_B$ and the polynomial degree (or moment order) $N_M$ are selected by scanning over a predefined grid $(N_B, N_M)$ and minimizing the two-sample Kolmogorov--Smirnov (K--S) test statistic between the empirical and reconstructed cumulative distribution functions. Performance of the proposed quantile-based piecewise monomial and Lagrange polynomial approaches are demonstrated on standard benchmark distributions (Normal, Weibull, bimodal Normal, trimodal Weibull) and on two empirical data sets, which are household electricity consumption and solar irradiance (via clear-sky index). Comparisons against KDE, SDE, standard monomial proposed in \citep{munkhammar2017polynomial} and Lagrange polynomial approaches proposed in \cite{turan2024polynomial} show that the quantile-based piecewise monomial and Lagrange polynomial approaches can achieve competitive or superior goodness-of-fit with reduced oscillations and improved robustness to skewness and multimodality.

\section{Moment-based piecewise polynomial reconstruction via quantile binning}

This section presents the general framework for constructing moment-based piecewise polynomial density estimations using quantile-based binning. The binning and a piecewise polynomial density estimation are described first, followed by the two polynomial formulations: piecewise monomial and piecewise Lagrange polynomials.

\subsection{Quantile-based binning and piecewise polynomial density estimation}

The observed data sample is denoted by $X = \{x_1, x_2, \dots, x_N\}$, and the aim is to estimate an underlying probability density function (PDF) $f(x)$ on a compact interval covering its support. The first step is to partition this support into $N_B$ bins. Define a sequence of bin boundaries $\mathbf{x}_{(1)} < \mathbf{x}_{(2)} < \dots < \mathbf{x}_{(N_B)}$,
chosen, in practice, as empirical quantiles such as
\begin{equation}
\mathbf{x}_{(i)} \approx F_{\mathrm{emp}}^{-1}\!\left(\frac{i}{N_B}\right), \qquad i = 1,\dots, N_B,    
\end{equation}
where $F_{\mathrm{emp}}$ is the empirical cumulative distribution function of the sample. This quantile-based construction yields bins containing approximately $N/N_B$ observations (up to ties and rounding).

Therefore, the disjoint bins $\{\Xi_i\}_{i=1}^{N_B}$ are defined by
\begin{equation} \label{bin_range}
    \Xi_i =
    \begin{cases}
    \{x: x \le \mathbf{x}_{(1)}\}, & i = 1, \\
    \{x: \mathbf{x}_{(i-1)} < x \le \mathbf{x}_{(i)}\}, & 1 < i < N_B, \\
    \{x: x > \mathbf{x}_{(N_B-1)}\}, & i = N_B
    \end{cases}
\end{equation}
such that $|\Xi_i|$ denotes the number of sample points falling in bin $\Xi_i$, and  $a_i = \min(\Xi_i)$ and $b_i = \max(\Xi_i)$ represent the lower and upper data-based bounds of the bin.

For each bin $\Xi_i$, empirical raw moments up to order $N_M - 1$ are computed:
\begin{equation}
    \tilde{m}_{i,j}
    = \frac{1}{|\Xi_i|} \sum_{x_k \in \Xi_i} x_k^{\,j}, 
    \qquad j = 0,1,\dots, N_M-1, \quad i = 1,2,\dots,N_B.
\end{equation}
The corresponding moment vector in the bin $\Xi_i$ is
\begin{equation}\label{dismom}
     \tilde{\mathbf{E}}_i = 
     \begin{bmatrix}
        \tilde{m}_{i,0} \\
        \tilde{m}_{i,1} \\
        \vdots \\
        \tilde{m}_{i,j}\\
        \vdots\\
        \tilde{m}_{i,N_M-1}
     \end{bmatrix},
\end{equation}
where $j$ is the order of moments. The density $f(x)$ is then approximated by a piecewise polynomial function
\begin{equation}
f(x) =
\begin{cases}
f_1(x), &  x \le \mathbf{x}_{(1)}, \\
f_i(x), & \mathbf{x}_{(i-1)} < x \le \mathbf{x}_{(i)}, \quad i = 2,\dots,N_B-1,\\
f_{N_B}(x), &  x > \mathbf{x}_{(N_B-1)}.
\end{cases}
\end{equation}
Each piecewise density $f_i(x)$ is supported (in practice) on $[a_i,b_i]$ and is approximated within that interval by a polynomial of degree at $N_M-1$. The coefficients of the piecewise polynomials are determined by matching the first $N_M$ empirical moments in each bin equalizing the empirical and approximated moments. This construction yields a piecewise mixture representation of $f(x)$ in which each bin contributes a local polynomial component based on balanced sample sizes.

Because the bins are defined via quantiles, they contain approximately equal numbers of observations. This reduces the risk of severe data sparsity in the tails and gives each bin a comparable amount of information for estimating local moments. In the case $N_B=2$, the formulation reduces to the standard monomial basis \cite{munkhammar2017polynomial}. When Lagrange polynomials constructed on Chebyshev nodes are used instead of monomials, the formulation corresponds to the Lagrange polynomial approach \cite{turan2024polynomial}.

\subsection{Piecewise monomial formulation}

In the piecewise monomial formulation, the local density in bin $\Xi_i$ is represented by a polynomial in the monomial basis \cite{munkhammar2017polynomial}:
\begin{equation}
  f_i(x) \approx P_i(x) = w_{i,0} + w_{i,1} x + \dots + w_{i,N_M-1} x^{N_M-1},
  \qquad x \in [a_i,b_i],
\end{equation}
where $\mathbf{W}_i = (w_{i,0},w_{i,1},\dots,w_{i,N_M-1})^\top$ is the vector of unknown coefficients for bin $i$.

The approximated moments of order $j$ in bin $\Xi_i$ with respect to $P_i(x)$ are
\begin{equation}
   \hat{m}_{i,j}
   = \int_{a_i}^{b_i} x^{j} P_i(x)\,dx
   = \int_{a_i}^{b_i} x^j
   \bigl(w_{i,0} + w_{i,1}x + \dots + w_{i,N_M-1}x^{N_M-1}\bigr)\,dx,
   \quad j = 0,\dots,N_M-1.
\end{equation}
By expanding the integrand, this can be written as a linear system in $\mathbf{W}_i$:
\begin{equation} \label{pwmonosystem}
    \mathbf{M}_i \mathbf{W}_i = \mathbf{E}_i, \quad i=1,2,\dots,N_B,
\end{equation}
where $\mathbf{E}_i$ is the vector of approximated moments (which are matched to empirical moments) and the moment matrix $\mathbf{M}_i$ has entries corresponding to integrals of monomials on $[a_i,b_i]$. Specifically,
\begin{equation}
    \mathbf{M}_i =
    \begin{bmatrix}
        b_i-a_i & \dfrac{b_i^2-a_i^2}{2} & \ldots & \dfrac{b_i^{N_M}-a_i^{N_M}}{N_M} \\
        \dfrac{b_i^2-a_i^2}{2} & \dfrac{b_i^3-a_i^3}{3} & \ldots & \dfrac{b_i^{N_M+1}-a_i^{N_M+1}}{N_M+1} \\
        \vdots & \vdots & \ddots & \vdots \\
        \dfrac{b_i^{N_M}-a_i^{N_M}}{N_M} &
        \dfrac{b_i^{N_M+1}-a_i^{N_M+1}}{N_M+1} &
        \ldots &
        \dfrac{b_i^{2N_M-1}-a_i^{2N_M-1}}{2N_M-1}
    \end{bmatrix},
\end{equation}
and
\begin{equation}
    \mathbf{E}_i =
    \begin{bmatrix}
        E_{\Xi_i}[x^0] \\
        E_{\Xi_i}[x^1] \\
        \vdots \\
        E_{\Xi_i}[x^{N_M-1}]
    \end{bmatrix}
    =
    \tilde{\mathbf{E}}_i,
\end{equation}
i.e., approximated moments are equated to empirical moments in bin $\Xi_i$ as in \eqref{dismom}. Solving \eqref{pwmonosystem} for each bin yields the local coefficient vectors $\mathbf{W}_i$ that define the piecewise monomial approximation $P(x)$ over the entire support. The procedure of piecewise monomial approach is described in Algorithm~\ref{monoalg}.

\begin{algorithm}[H]
\caption{Piecewise monomial procedure} \label{monoalg}
\begin{algorithmic}[1]
\State \textbf{Input:} Data $X = \{x_1,\dots,x_N\}$
\State Choose a goodness-of-fit test $\mathbf{T}$ (e.g., K--S test statistic)
\State Define range of bins $N_{B,\text{list}} = \{1,2,\dots,19\}$ and range of number of moments $N_{M,\text{list}} = \{3,4,\dots,11\}$
\For{each $N_B \in N_{B,\text{list}}$}
    \State Compute empirical quantile bin boundaries and form bins $\Xi_i$, $i=1,\dots,N_B$ \eqref{bin_range}
    \State Extract $a_i = \min(\Xi_i)$ and $b_i = \max(\Xi_i)$ for each bin
    \For{each $N_M \in N_{M,\text{list}}$}
        \State Compute empirical moment vectors $\tilde{\mathbf{E}}_i$ in each bin using $\tilde{m}_{i,j}$
        \State Assemble $\mathbf{M}_i$ from $a_i,b_i$ and solve $\mathbf{M}_i \mathbf{W}_i = \tilde{\mathbf{E}}_i$
        \State Construct the piecewise polynomial $P(x)$ using $\{\mathbf{W}_i\}_{i=1}^{N_B}$ \eqref{pwmonosystem}
        \State Compute the empirical CDF $F_{\text{emp}}$ and the CDF $F_P$ of $P(x)$
        \State Evaluate $\mathbf{T}(N_B, N_M) = \sup_x |F_P(x) - F_{\text{emp}}(x)|$
    \EndFor
\EndFor
\State Identify feasible combinations where $P(x) \ge 0$ (or nearly non-negative) on the support
\State Among feasible $(N_B,N_M)$, choose $(N_B^*, N_M^*)$ that minimize $\mathbf{T}(N_B,N_M)$
\State \textbf{Output:} Selected $(N_B^*, N_M^*, \mathbf{T}^*)$ and corresponding $P^*(x)$
\end{algorithmic}
\end{algorithm}
In the procedure described in Algorithm~\ref{monoalg}, the feasibility condition $P(x) \ge 0$ ensures that the reconstructed function is a valid density (or nearly so, up to a valid density). If strict non-negativity cannot be achieved across the grid, the combination that is closest to non-negative while minimizing the goodness-of-fit statistic is selected.

\subsection{Piecewise Lagrange polynomial formulation}

As an alternative to the monomial basis, a piecewise Lagrange polynomial representation is considered. In bin $\Xi_i$, the local density is approximated by \cite{turan2024polynomial}
\begin{equation}
  f_i(x) \approx L_i(x)
  = \omega_{i,0} \, l_{i,0}(x)
  + \omega_{i,1} \, l_{i,1}(x)
  + \dots
  + \omega_{i,N_M-1} \, l_{i,N_M-1}(x),
  \qquad x \in [a_i,b_i],
\end{equation}
where $\hat{\mathbf{W}}_i = (\omega_{i,0},\dots,\omega_{i,N_M-1})^\top$ are unknown coefficients and $l_{i,k}(x)$ are Lagrange basis polynomials associated with interpolation points $\{x_{i,\lambda}\}_{\lambda=0}^{N_M-1}$ in bin $i$.

The interpolation points are chosen as Chebyshev nodes mapped to $[a_i,b_i]$:
\begin{equation}
    x_{i,\lambda}
    = a_i + \frac{b_i - a_i}{2}
      \left[
        1 - \cos\!\left(\frac{\lambda \pi}{N_M - 1}\right)
      \right],
    \quad \lambda = 0,1,\dots,N_M-1.
\end{equation}
For each bin $\Xi_i$, the Lagrange basis polynomials are then defined by
\begin{equation}
    l_{i,k}(x)
    = \prod_{\lambda = 0,\ \lambda \ne k}^{N_M-1}
      \frac{x - x_{i,\lambda}}{x_{i,k} - x_{i,\lambda}},
    \quad k = 0,\dots,N_M-1.
\end{equation}

The approximated moments in bin $\Xi_i$ induced by $L_i(x)$ are
\begin{equation}
    \hat{E}_{\Xi_i}[x^j]
    = \int_{a_i}^{b_i} x^{j} L_i(x)\,dx
    = \int_{a_i}^{b_i} x^{j}
      \Bigl(
        \omega_{i,0} l_{i,0}(x) + \dots + \omega_{i,N_M-1} l_{i,N_M-1}(x)
      \Bigr)\,dx,
    \quad j = 0,\dots,N_M-1.
\end{equation}
These integrals define the entries of a linear system
\begin{equation}\label{eqsystemLag}
   \mathbf{L}_i \hat{\mathbf{W}}_i = \hat{\mathbf{E}}_i, \quad i = 1,2,\dots,N_B,
\end{equation}
where
\begin{equation}
   \mathbf{L}_i =
   \begin{bmatrix}
       \displaystyle \int_{a_i}^{b_i} l_{i,0}(x)\,dx
         & \displaystyle \int_{a_i}^{b_i} l_{i,1}(x)\,dx
         & \ldots
         & \displaystyle \int_{a_i}^{b_i} l_{i,N_M-1}(x)\,dx \\
       \displaystyle \int_{a_i}^{b_i} x\,l_{i,0}(x)\,dx
         & \displaystyle \int_{a_i}^{b_i} x\,l_{i,1}(x)\,dx
         & \ldots
         & \displaystyle \int_{a_i}^{b_i} x\,l_{i,N_M-1}(x)\,dx \\
       \vdots & \vdots & \ddots & \vdots \\
       \displaystyle \int_{a_i}^{b_i} x^{N_M-1} l_{i,0}(x)\,dx
         & \displaystyle \int_{a_i}^{b_i} x^{N_M-1} l_{i,1}(x)\,dx
         & \ldots
         & \displaystyle \int_{a_i}^{b_i} x^{N_M-1} l_{i,N_M-1}(x)\,dx
   \end{bmatrix},
\end{equation}
and
\begin{equation}
    \hat{\mathbf{E}}_i =
    \begin{bmatrix}
        \hat{E}_{\Xi_i}[x^0] \\
        \hat{E}_{\Xi_i}[x^1] \\
        \vdots \\
        \hat{E}_{\Xi_i}[x^{N_M-1}]
    \end{bmatrix}
    =
    \tilde{\mathbf{E}}_i.
\end{equation}
As before, approximated moments are equated to empirical moments bin-wise. Solving \eqref{eqsystemLag} yields $\hat{\mathbf{W}}_i$ and thus the local Lagrange polynomial approximation $L_i(x)$ in each bin. The procedure of piecewise Lagrange polynomial approach is described in Algorithm~\ref{Lagrangealg}.

\begin{algorithm}[H]
\caption{Piecewise Lagrange polynomial procedure} \label{Lagrangealg}
\begin{algorithmic}[1]
\State \textbf{Input:} Data $X = \{x_1,\dots,x_N\}$
\State Choose a goodness-of-fit test $\mathbf{T}$ (e.g., K--S test statistic)
\State Define range of bins $N_{B,\text{list}} = \{1,2,\dots,19\}$ and range of number of moments $N_{M,\text{list}} = \{3,4,\dots,11\}$
\For{each $N_B \in N_{B,\text{list}}$}
    \State Compute empirical quantile bin boundaries and form bins $\Xi_i$, $i=1,\dots,N_B$
    \State Extract $a_i = \min(\Xi_i)$ and $b_i = \max(\Xi_i)$ for each bin
    \For{each $N_M \in N_{M,\text{list}}$}
        \State Compute empirical moment vectors $\tilde{\mathbf{E}}_i$ in each bin
        \State Construct Chebyshev interpolation points $\{x_{i,\lambda}\}_{\lambda=0}^{N_M-1}$ on $[a_i,b_i]$
        \State Build Lagrange basis polynomials $\{l_{i,k}(x)\}_{k=0}^{N_M-1}$ and form matrix $\mathbf{L}_i$
        \State Solve $\mathbf{L}_i \hat{\mathbf{W}}_i = \tilde{\mathbf{E}}_i$ for each bin
        \State Construct the piecewise polynomial $L(x)$ from $\{\hat{\mathbf{W}}_i\}_{i=1}^{N_B}$
        \State Compute the empirical CDF $F_{\text{emp}}$ and the CDF $F_L$ of $L(x)$
        \State Evaluate $\mathbf{T}(N_B, N_M) = \sup_x |F_L(x) - F_{\text{emp}}(x)|$
    \EndFor
\EndFor
\State Identify feasible combinations where $L(x) \ge 0$ (or nearly non-negative) on the support
\State Among feasible $(N_B,N_M)$, choose $(N_B^*, N_M^*)$ that minimize $\mathbf{T}(N_B,N_M)$
\State \textbf{Output:} Selected $(N_B^*, N_M^*, \mathbf{T}^*)$ and corresponding $L^*(x)$
\end{algorithmic}
\end{algorithm}

\section{Applications}\label{sec-verify}

This section evaluates the proposed piecewise polynomial procedures with quantile-based binning on benchmark distributions and real-world data sets. The approaches, which are Kernel Density Estimation (KDE), Spline Density Estimation (SDE), standard monomial \citep{munkhammar2017polynomial}, Lagrange polynomial \citep{turan2024polynomial}, piecewise monomial (proposed), and piecewise Lagrange polynomial (proposed) are compared.

For quantitative assessment, the two-sample Kolmogorov--Smirnov (K--S) test statistic \cite{conover1999practical} is used:
\begin{equation}
    K\!-\!S = \sup_{\xi} \bigl| F_{\text{emp}}(\xi) - F_{P}(\xi) \bigr|,
\end{equation}
where $F_{P}$ denotes the cumulative density function (CDF) associated with a given reconstruction method (here $P$ stands for a reconstructed PDF). In addition, a simple goodness-of-fit index \cite{biggs2002modelling} is considered:
\begin{equation}
    GoF = \frac{\hat{F}_P - st_{\text{error}}}{\hat{F}_P},
\end{equation}
where $\hat{F}_P$ is the mean of the estimated distribution $F_P$ and
\begin{equation}
   st_{\text{error}} =
   \sqrt{
     \frac{1}{N - N_B N_M}
     \sum_{i=1}^{N}
       \bigl(F_{\text{emp}}(x_i) - F_{P}(x_i)\bigr)^2
   }.
\end{equation}
The $GoF$ index is closer to 1 when the estimator matches the empirical CDF more closely (small $st_{\text{error}}$), and it degrades as the discrepancy increases.

\subsection{Estimation of well-known distribution functions}

The first test concerns synthetic data drawn from four benchmark distributions: Normal, Weibull, bimodal Normal, and trimodal Weibull. For each distribution, $N$ samples are generated and the true density is known in closed form.

\begin{figure}[H]
    \centering
    \includegraphics[width=1.1\linewidth]{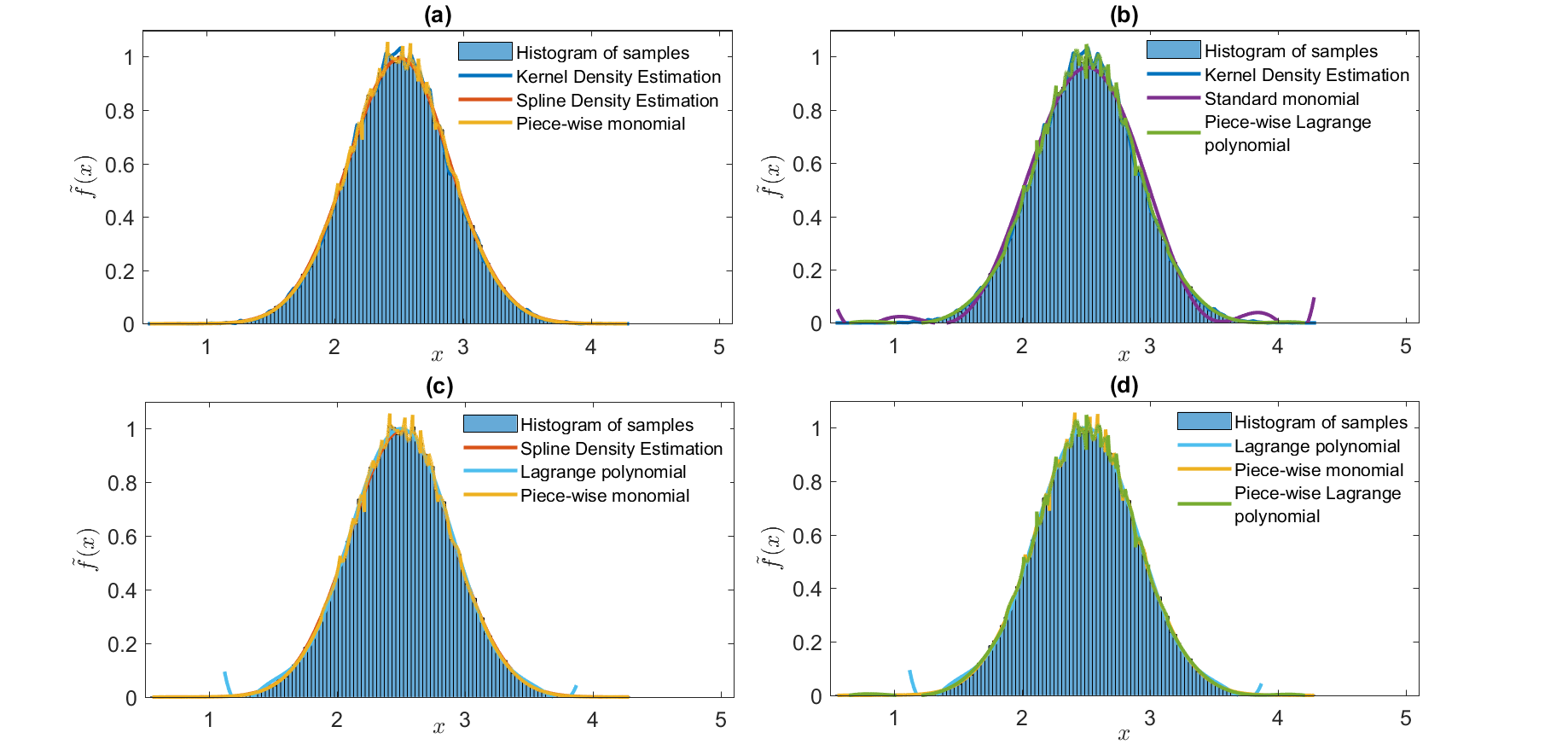}
    \caption{Normal distribution: Bandwidth $h=0.05$ for KDE; $NumBasis=14$ and $Degree=3$ for SDE; $N_M=11$ for standard monomial and Lagrange polnomial approaches; $N_M=5$, $N_B=18$ for piecewise monomial; and $N_M=4$, $N_B=14$ for piecewise Lagrange polynomial.}
    \label{fig:Normal_Dist}
\end{figure}

\begin{figure}[H]
    \centering
    \includegraphics[width=1.05\linewidth]{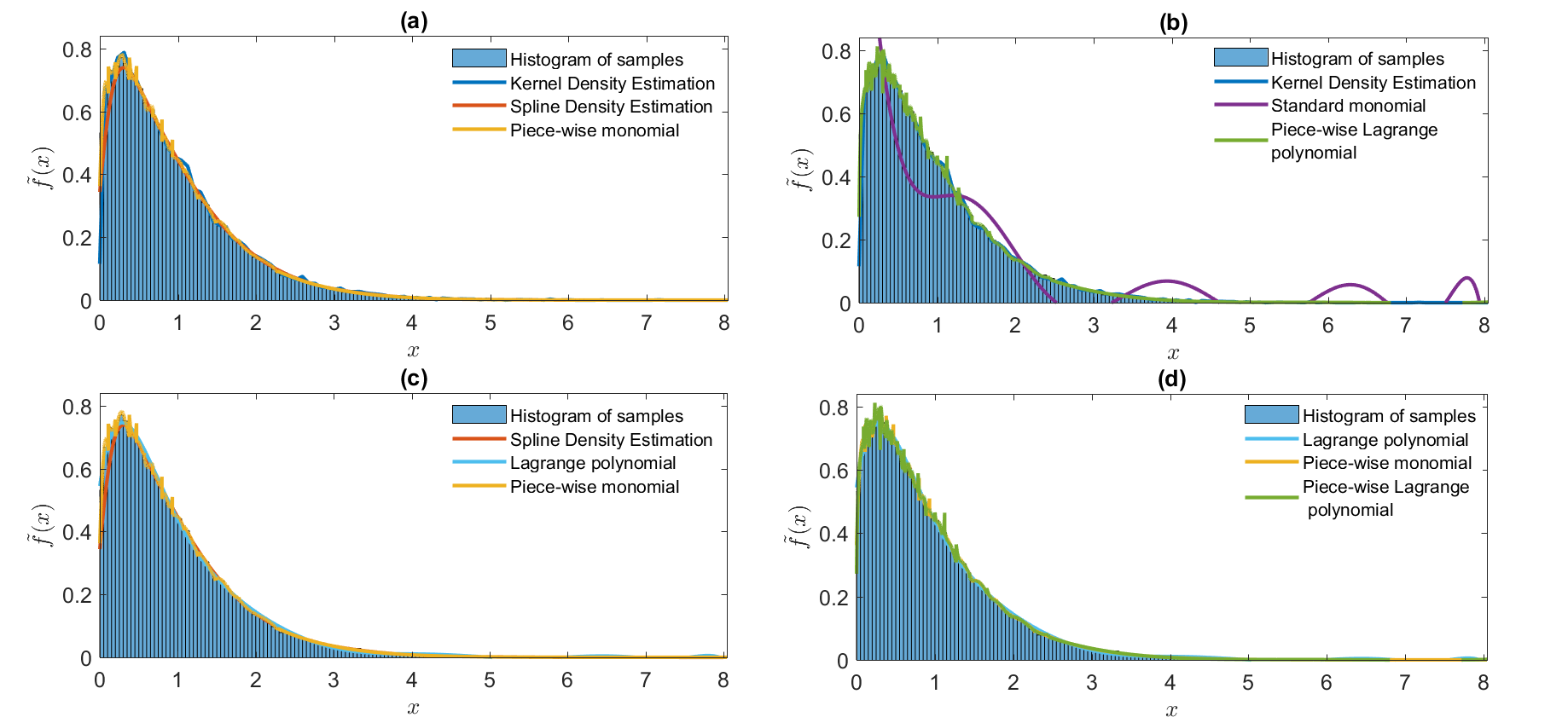}
    \caption{Weibull distribution: Bandwidth $h=0.05$ for KDE; $NumBasis=20$ and $Degree=3$ for SDE; $N_M=11$ for standard monomial and Lagrange polynomial approaches; $N_M=7$, $N_B=19$ for piecewise monomial; and $N_M=5$, $N_B=16$ for piecewise Lagrange polynomial.}
    \label{fig:Weibull}
\end{figure}

\begin{figure}[H]
    \centering
    \includegraphics[width=1\linewidth]{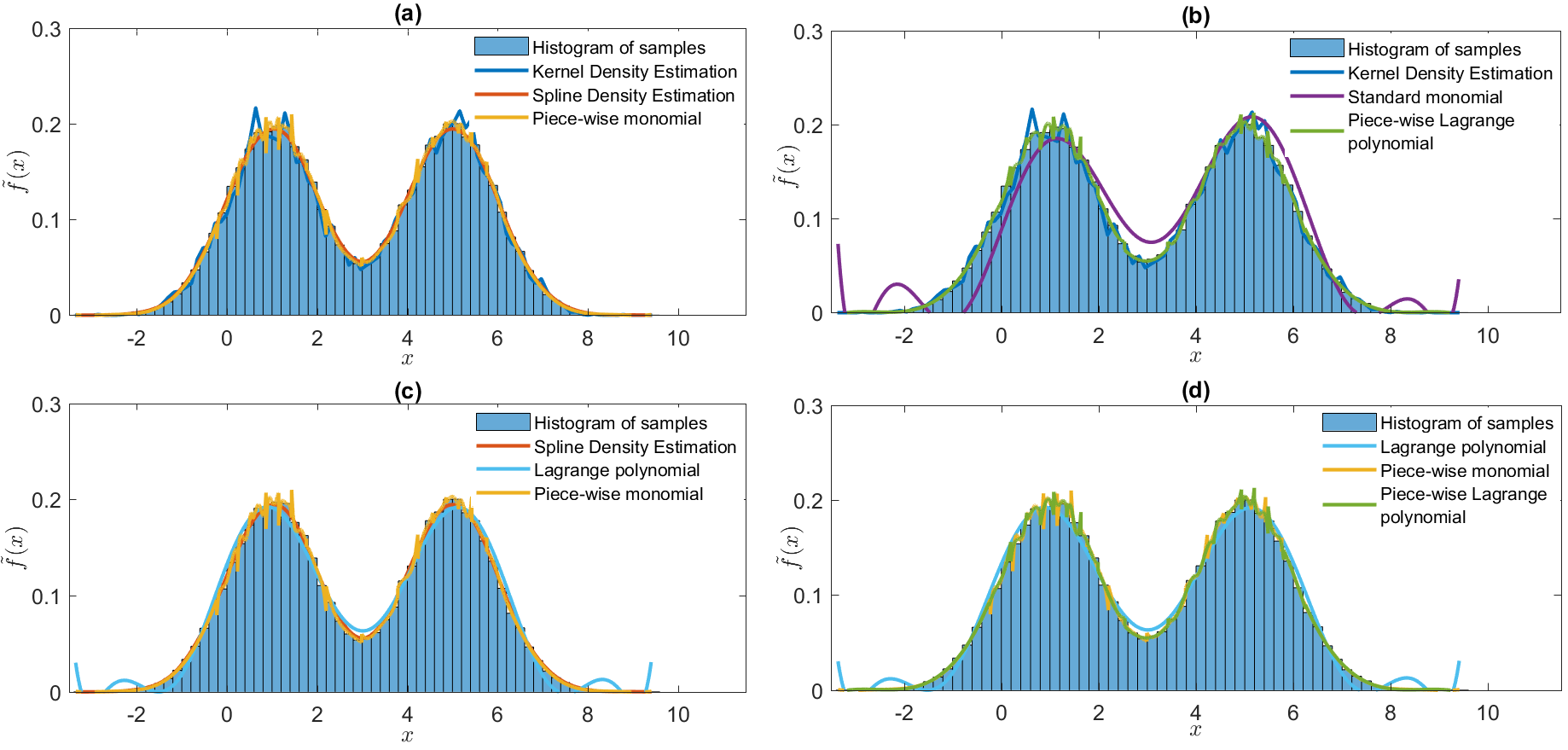}
    \caption{Bimodal Normal distribution: Bandwidth $h=0.05$ for KDE; $NumBasis=19$ and $Degree=3$ for SDE; $N_M=11$ for standard monomial and Lagrange polynomial approaches; $N_M=8$, $N_B=16$ for piecewise monomial; and $N_M=5$, $N_B=14$ for piecewise Lagrange polynomial.}
    \label{fig:Bimodal_Normal}
\end{figure}

\begin{figure}[H]
    \centering
    \includegraphics[width=1\linewidth]{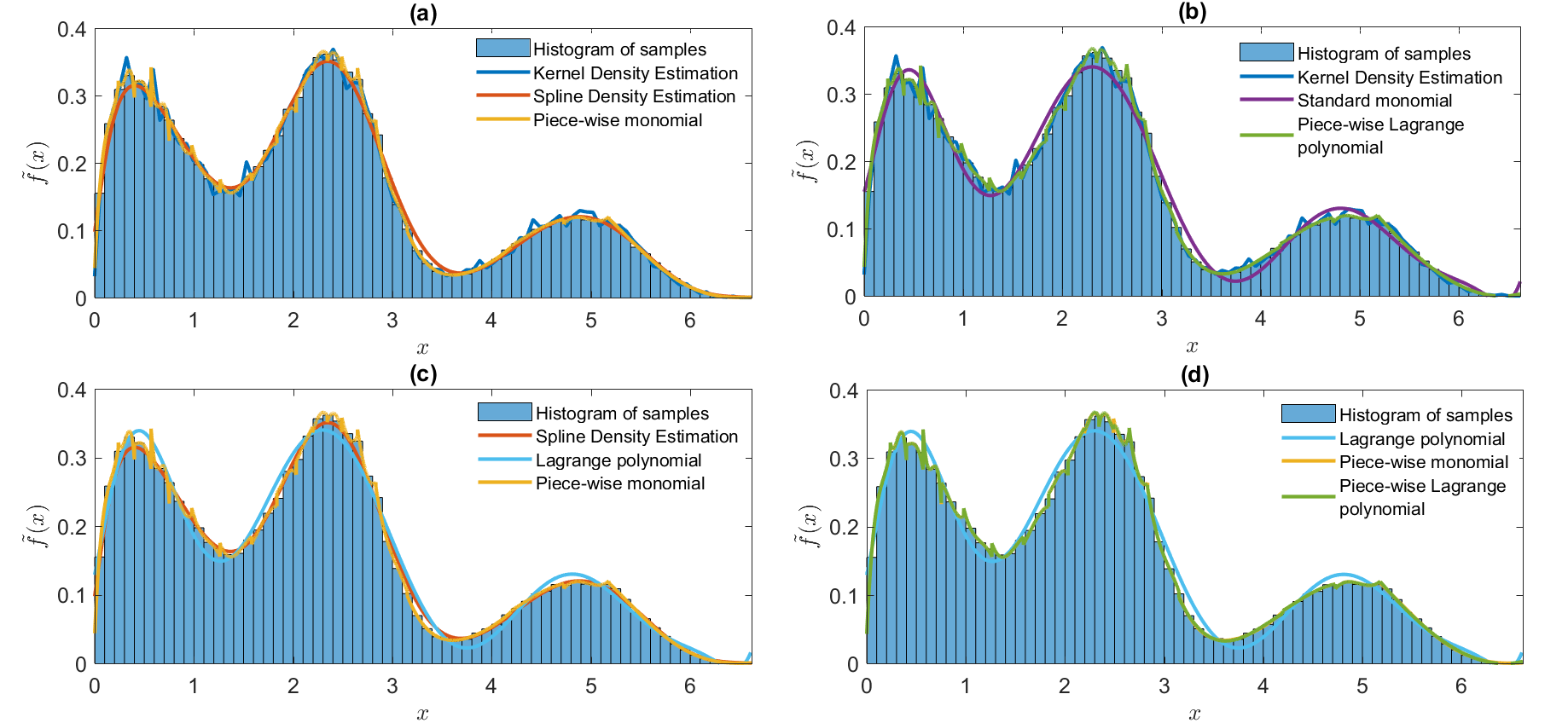}
    \caption{Trimodal Weibull distribution: Bandwidth $h=0.05$ for KDE; $NumBasis=18$ and $Degree=7$ for SDE; $N_M=11$ for standard monomial and Lagrange polnomial approaches; $N_M=5$, $N_B=19$ for piecewise monomial; and $N_M=4$, $N_B=19$ for piecewise Lagrange polynomial.}
    \label{fig:Trimodal_Weibull}
\end{figure}

Figures~\ref{fig:Normal_Dist}--\ref{fig:Trimodal_Weibull} show the reconstructed PDFs for all six methods:
Kernel Density Estimation (KDE), Spline Density Estimation (SDE), standard monomial \cite{munkhammar2017polynomial}, Lagrange polynomial \cite{turan2024polynomial}, piecewise monomial, and piecewise Lagrange polynomial. For each case, the bandwidth used in KDE, the number of basis functions and degree used in SDE, and the selected $(N_B,N_M)$ for the piecewise approaches are reported.

As shown in Figures~\ref{fig:Normal_Dist}–\ref{fig:Trimodal_Weibull}, the piecewise approaches accurately capture both the peak regions and the distribution tails, whereas the monomial and Lagrange polynomial approaches exhibit pronounced oscillations in the tail regions.

Table~\ref{KSwellknown} summarizes the K--S test statistics for all approaches and all four benchmark distributions, together with the distribution parameters used in this study. Lower K--S test statistic values indicate better agreement between the empirical and reconstructed CDFs.

\begin{table}[H]
\caption{$K\!-\!S$ test statistics of KDE, SDE, standard monomial, Lagrange polynomial, piecewise monomial, and piecewise Lagrange polynomial approaches for well-known distributions.} \label{KSwellknown}
\centering
\scriptsize
\begin{tabular}{c|c|c|c|c|c|c|c}
       & parameters & \makecell{KDE} & SDE &
       \makecell{standard\\monomial} &
       \makecell{Lagrange\\polynomial} &
       \makecell{piecewise\\monomial} &
       \makecell{piecewise Lagrange\\ polynomial} \\
\hline
Normal & $\sigma /\mu=0.16$ &4.12e-04 &2.90e-03 &1.03e-02 &8.72e-03 &5.52e-04 & 7.34e-04\\ \hline 
Weibull &$k=1.0$, $\lambda=1.2$ &7.35e-04 &1.80e-02 &3.72e-03 &2.82e-03 &5.45e-04 &4.79e-04 \\ \hline
\makecell{Bimodal\\Normal} & \makecell{$\sigma_1 /\mu_1=1.0$\\ $\sigma_2 /\mu_2=0.2$} &1.83e-04 &2.60e-03 &9.65e-03 &2.14e-03 &4.33e-04 &4.83e-04 \\  \hline
\makecell{Trimodal\\Weibull} &\makecell{$k=0.8$, $\lambda=1.5$ \\ $k=2.5$, $\lambda=5.2$ \\ $k=5.0$, $\lambda=8.2$} & 2.70e-04&8.03e-03 &7.08e-03 &7.50e-03 &4.90e-04 & 4.34e-04\\  
    \end{tabular}
\end{table}

\begin{table}[H]
\centering
\small
\caption{Percentage improvement of the piecewise monomial approach relative to other approaches.}
\label{tab:pwmon_vs_others}
\begin{tabular}{c|rrrrr}
\hline
Distribution & KDE (\%)& SDE (\%) & \makecell{standard\\monomial (\%)} & \makecell{Lagrange\\polynomial (\%)}&
\makecell{piecewise Lagrange\\polynomial (\%)}\\
\hline
Normal           & -33.98\% & 80.97\% & 94.64\% & 93.67\% & 24.80\% \\
Weibull          &  25.85\% & 96.97\% & 85.35\% & 80.67\% & -13.78\% \\
Bimodal Normal   & -136.61\% & 83.35\% & 95.51\% & 79.77\% & 10.35\% \\
Trimodal Weibull & -81.48\% & 93.90\% & 93.08\% & 93.47\% & -12.90\% \\
\hline
\end{tabular}
\end{table}

\begin{table}[H]
\centering
\small
\caption{Percentage improvement of the piecewise Lagrange polynomial approach relative to other approaches.}
\label{tab:pwlag_vs_others}
\begin{tabular}{c|rrrrr}
\hline
Distribution & KDE (\%) & SDE (\%) & \makecell{standard\\monomial (\%)} & \makecell{Lagrange\\polynomial (\%)}&
\makecell{piecewise\\monomial (\%)}\\
\hline
Normal           & -78.16\% & 74.69\% & 92.87\% & 91.58\% & -32.97\% \\
Weibull          &  34.83\% & 97.34\% & 87.12\% & 83.01\% & 12.11\% \\
Bimodal Normal   & -163.93\% & 81.42\% & 94.99\% & 77.43\% & -11.55\% \\
Trimodal Weibull & -60.74\% & 94.60\% & 93.87\% & 94.21\% & 11.43\% \\
\hline
\end{tabular}
\end{table}

Tables \ref{tab:pwmon_vs_others} and \ref{tab:pwlag_vs_others} summarize the relative improvements in K–S statistics obtained by the two proposed approaches—piecewise monomial and piecewise Lagrange polynomial—compared with the benchmark density–estimation techniques for four well-known distributions with increasing multimodality. For all test cases, the piecewise monomial representation yields an improvement of 80–96\% relative to standard monomial and Lagrange polynomial approaches. Similarly, the piecewise Lagrange polynomial approach shows 77–95\% improvement over standard monomial and Lagrange polynomial approaches. These results confirm that localizing the approximation is crucial for reducing oscillations and improving the fidelity of polynomial-based PDF reconstruction, especially in multimodal distributions where global polynomials struggle. Second, both piecewise approaches consistently outperform SDE by large margins. Improvements exceed 75–97\%, demonstrating that SDE remains less accurate for sharply varying or multimodal PDFs, where spline smoothness constraints tend to oversmooth the density and distort local extrema. Third, comparisons against KDE reveal that the piecewise approaches do not always dominate. KDE exhibits very small K–S errors for the Normal and Bimodal Normal cases, leading to negative improvement percentages for both piecewise approaches in these instances. This reflects the well-known variance–bias trade-off in KDE: when the underlying distribution is smooth and well-resolved by data, KDE may achieve exceptional accuracy. However, in the Weibull case, both piecewise approaches outperform KDE by 25–35\%, illustrating their advantage in handling skewed shape. When directly comparing the two proposed methods, piecewise monomial and piecewise Lagrange polynomial achieve comparable performance. The piecewise Lagrange polynomial approach is moderately more accurate in most cases (approximately 12\% improvement over piecewise monomial for the Weibull and Trimodal Weibull distributions), whereas piecewise monomial slightly outperform piecewise Lagrange polynomial in the Normal and Bimodal Normal cases (10-25\%). The reason is the superiority of Lagrange polynomials with Chebyshev nodes in reducing oscillations and capturing skewed shape. These differences, however, remain relatively small compared with the gains relative to standard monomial, Lagrange polynomial and SDE approaches. 

\subsection{Estimation of empirical distribution functions}

Two real-world data sets are now considered to illustrate the practical relevance of the proposed approach.

\paragraph{Household electricity consumption.}
The first data set consists of measured household electricity consumption with a ten-minute temporal resolution for a detached house over one year \cite{munkhammar2014characterizing}. This type of data is often modeled using unimodal parametric families such as Log-Normal or Weibull distributions. However, real household electricity consumption profiles typically deviate from simple unimodal shapes due to daily cycles, occupancy patterns, and appliance usage, leading to asymmetries and possible secondary modes.

\begin{figure}[H]
    \centering
    \includegraphics[width=1\linewidth]{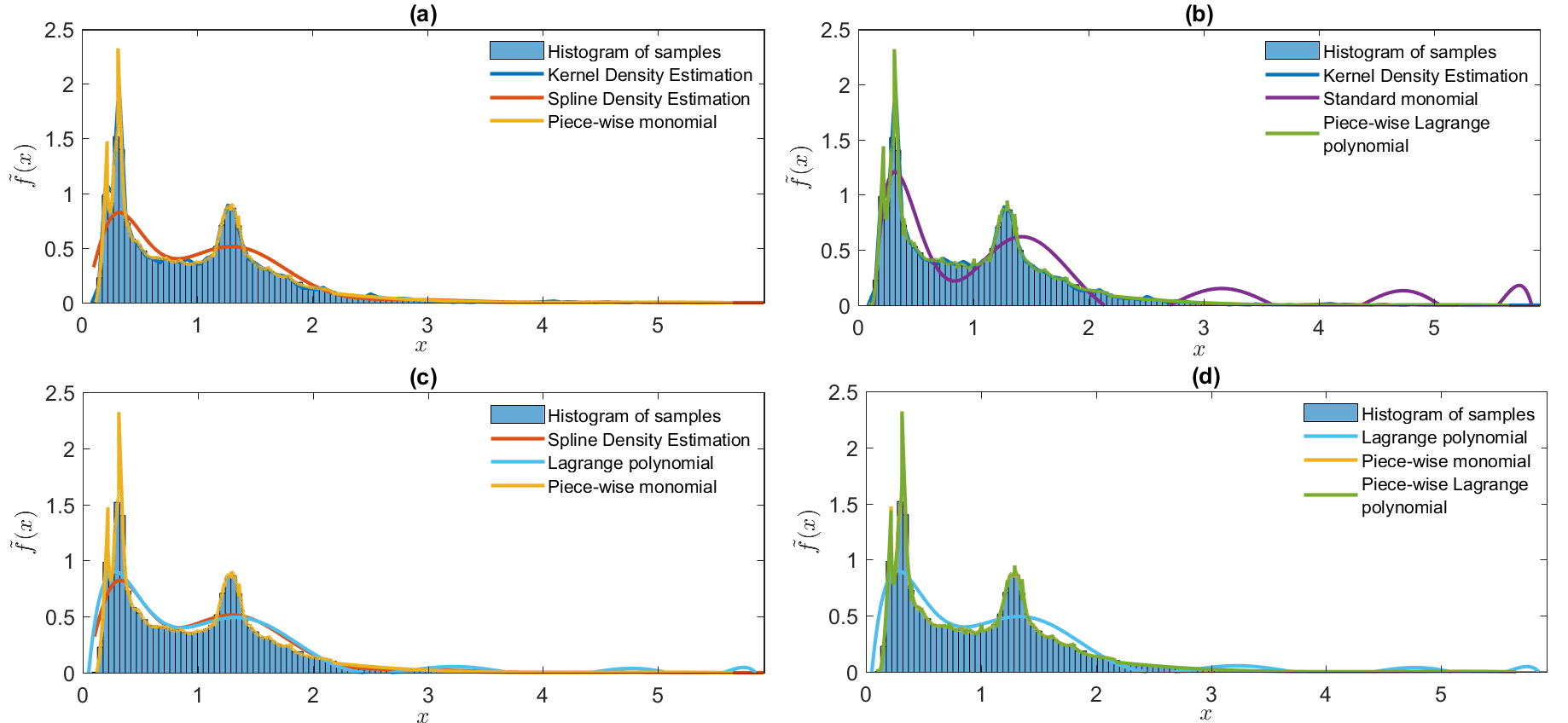}
    \caption{Household electricity consumption data set: Bandwidth $h=0.05$ for KDE; $NumBasis=17$ and $Degree=6$ for SDE; $N_M=11$ for standard monomial and Lagrange polynomial approaches; $N_M=4$, $N_B=19$ for piecewise monomial; and $N_M=4$, $N_B=19$ for piecewise Lagrange polynomial.}
    \label{fig:Electricity_Data}
\end{figure}

The household electricity consumption data are fitted using KDE, SDE, standard monomial \cite{munkhammar2017polynomial}, Lagrange polynomial \cite{turan2024polynomial}, and the two proposed piecewise approaches. Figure~\ref{fig:Electricity_Data} shows the resulting PDF estimates along with the empirical histogram. The caption lists the bandwidth for KDE, the number of spline basis functions and degree for SDE, and the selected $(N_B,N_M)$ combinations for the piecewise approaches.

The piecewise approaches adapt to local features of the household electricity consumption distribution and capture both the main body and tail behavior without introducing large oscillations. The polynomial methods, in contrast, show more pronounced wiggles near the boundaries. KDE and SDE provide smooth fits but require careful tuning of bandwidth or basis parameters.

\begin{table}[H]
\caption{K--S test statistics of different PDF reconstruction approaches for the household electricity consumption data set.}    \label{tab:ks_Elec}
    \centering
    \begin{tabular}{l c}
        \hline
        \textbf{Method} & \textbf{Household electricity consumption} \\
        \hline
        Kernel Density Estimation        & 1.08e-03  \\
        Spline Density Estimation        & 4.73e-02 \\
        standard monomial                  & 5.43e-02\\
        Lagrange polynomial       & 3.66e-02 \\
        piecewise monomial               & 1.75e-03  \\
        piecewise Lagrange polynomial    & 1.96e-03  \\
        \hline
    \end{tabular}
\end{table}

Table~\ref{tab:ks_Elec} reports the Kolmogorov–Smirnov (K-S) test statistics for different probability density function approaches applied to the household electricity consumption dataset. The results indicate that the piecewise approaches significantly outperform the monomial and Lagrange polynomial approaches, yielding K-S statistics comparable to kernel density estimation. Among the piecewise approaches, the piecewise monomial approach yields the lower K-S test statistic, suggesting a slightly improved fit relative to the piecewise Lagrange polynomial approach. 

\paragraph{Solar irradiance and clear-sky index.}
The second data set consists of radiometer measurements of global horizontal irradiance (GHI) for the year 2008, obtained from the Swedish Meteorological and Hydrological Institute (SMHI) for Norrköping, Sweden (59$^\circ$35$'$31$''$~N, 17$^\circ$11$'$8$''$~E) \cite{SMHI2008}. The original data provide instantaneous irradiance values with a one-minute resolution. To reduce the influence of low solar altitudes, 120 measurements (corresponding to 120 minutes) centered around noon are selected for each day, resulting in a total of 43681 observations. These irradiance values can be transformed to a clear-sky index (CSI) by suitable normalization, but here the focus is directly on the distribution of the selected GHI (or CSI) values \cite{munkhammar2018markov}.

\begin{figure}[H]
    \centering
    \includegraphics[width=1\linewidth]{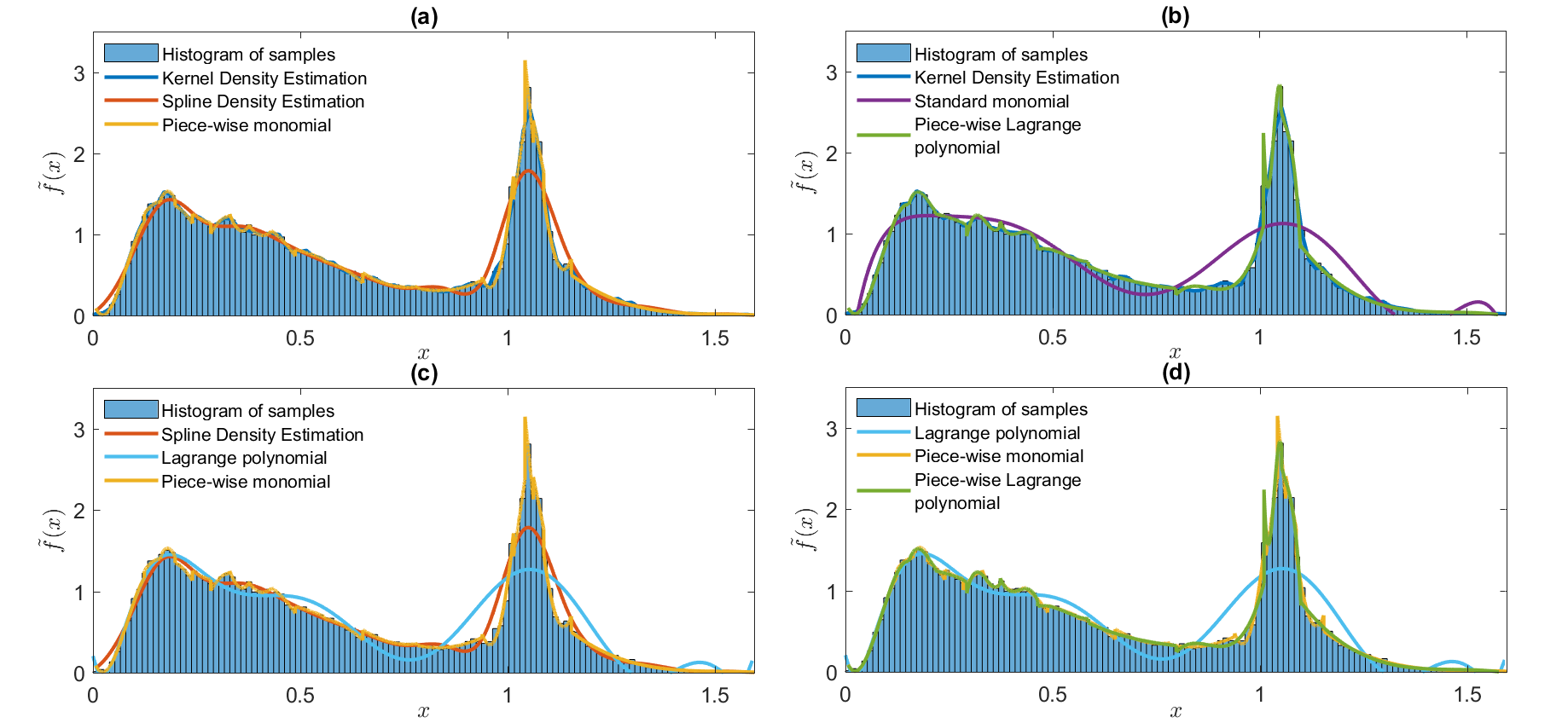}
    \caption{Solar irradiance (CSI) data set: Bandwidth $h=0.05$ for KDE; $NumBasis=20$ and $Degree=3$ for SDE; $N_M=11$ for standard monomial and Lagrange polynomial approaches; $N_M=4$, $N_B=19$ for piecewise monomial; and $N_M=4$, $N_B=11$ for piecewise Lagrange polynomial.}
    \label{fig:CSI_data}
\end{figure}

As for the household electricity data, KDE, SDE, standard monomial, Lagrange polynomial, and the proposed piecewise approaches are applied. Figure~\ref{fig:CSI_data} shows the reconstructed PDFs and the empirical histograms. Again, the caption specifies the hyperparameters for KDE and SDE and the chosen $(N_B,N_M)$ for the piecewise polynomial approaches.

Table~\ref{tab:ks_CSI} reports the K--S test statistics for all methods on solar irradiance data sets. The piecewise approaches achieve K--S test statistic values close to those of KDE, and significantly improve over polynomial approaches and SDE in these examples.

\begin{table}[H]
\caption{K--S test statistics of different PDF reconstruction approaches for solar irradiance (CSI) data set.}    \label{tab:ks_CSI}
    \centering
    \begin{tabular}{l cc}
        \hline
        \textbf{Method} & \textbf{Solar irradiance} \\
        \hline
        Kernel Density Estimation        & 2.16e-03 \\
        Spline Density Estimation         & 2.09e-02 \\
        standard monomial                   & 4.50e-02 \\
        Lagrange polynomial        & 3.84e-02 \\
        piecewise monomial                & 2.13e-03 \\
        piecewise Lagrange polynomial    & 1.92e-03 \\
        \hline
    \end{tabular}
\end{table}

\begin{table}[H]
\centering
\small
\caption{Percentage improvement of piecewise monomial and piecewise Lagrange polynomial relative to other PDF reconstruction approaches.} \label{data_per}
\begin{tabular}{l|cc|cc}
\hline
\multirow{2}{*}{\textbf{Method}} 
& \multicolumn{2}{c|}{\makecell{piecewise monomial\\(\%)}} 
& \multicolumn{2}{c}{\makecell{piecewise Lagrange\\polynomial (\%)}} \\
& \makecell{Household electricity\\consumption} & \makecell{Solar\\irradiance} & \makecell{Household electricity\\consumption} & \makecell{Solar\\irradiance}\\
\hline
Kernel Density Estimation        & -61.9\%  & 1.39\%  & $-81.9$\% & 11.1\% \\
Spline Density Estimation        & 96.3\%  & 89.8\%  & 95.9\% & 90.8\% \\
Standard monomial                & 96.8\%  & 95.3\%  & 96.4\% & 95.7\% \\
Lagrange polynomial              & 95.2\%  & 94.5\%  & 94.6\% & 95.0\% \\
\makecell{piecewise monomial\\(baseline)}   & ---     & ---     & $-12.0$\% & 9.9\% \\
\makecell{piecewise Lagrange\\polynomial (baseline)} & $10.7$\% & $-9.9$\% & --- & --- \\
\hline
\end{tabular}
\end{table}

Tables \ref{data_per} and the corresponding percentage improvement table summarize the performance of the PDF reconstruction approaches used in this study on household electricity consumption and solar irradiance data. Piecewise monomial and piecewise Lagrange polynomial approaches reduce the K--S test statistic by approximately 95\% for both data sets compared to SDE, standard monomial and Lagrange polynomial approaches. 

Compared to kernel density estimation (KDE), both piecewise approaches exhibit slightly inferior performance for the household electricity consumption dataset, while providing improved performance for the solar irradiance dataset. 

Between the two piecewise approaches, the piecewise monomial approach demonstrates superior performance for the household electricity consumption dataset, as evidenced by its lower K-S test statistic. This improvement can be attributed to its simpler local basis that mitigates oscillatory artifacts and enhances numerical robustness in smoothly varying regimes, considering the K-S test statistics. Household electricity consumption time series are commonly well described by a mixture of distributions belonging to the same statistical family, which favors stable low-order local approximations. In contrast, for the solar irradiance dataset, the piecewise Lagrange polynomial approach yields slightly better performance. This improvement is primarily due to the use of Chebyshev nodes, which minimize interpolation error and capture pick point. Such properties are especially advantageous for solar irradiance data, which are better represented as a mixture of fundamentally different distribution families, where the enhanced local flexibility of Lagrange polynomials constructed on Chebyshev nodes enables a more accurate representation of heterogeneous, asymmetric features compared to symmetric representations.

\subsection{Sensitivity analysis}

To better understand the behavior of the piecewise approaches, the dependence of the K--S test statistic on the number of bins $N_B$ and the polynomial moment order $N_M$ is investigated. For both empirical data sets, $N_B = 1,\dots,19$ and $N_M = 3,\dots,11$ are scanned and the K--S test statistic is computed for each combination, for both the piecewise approaches.

\begin{figure}[H]
    \centering
    \includegraphics[width=1\linewidth]{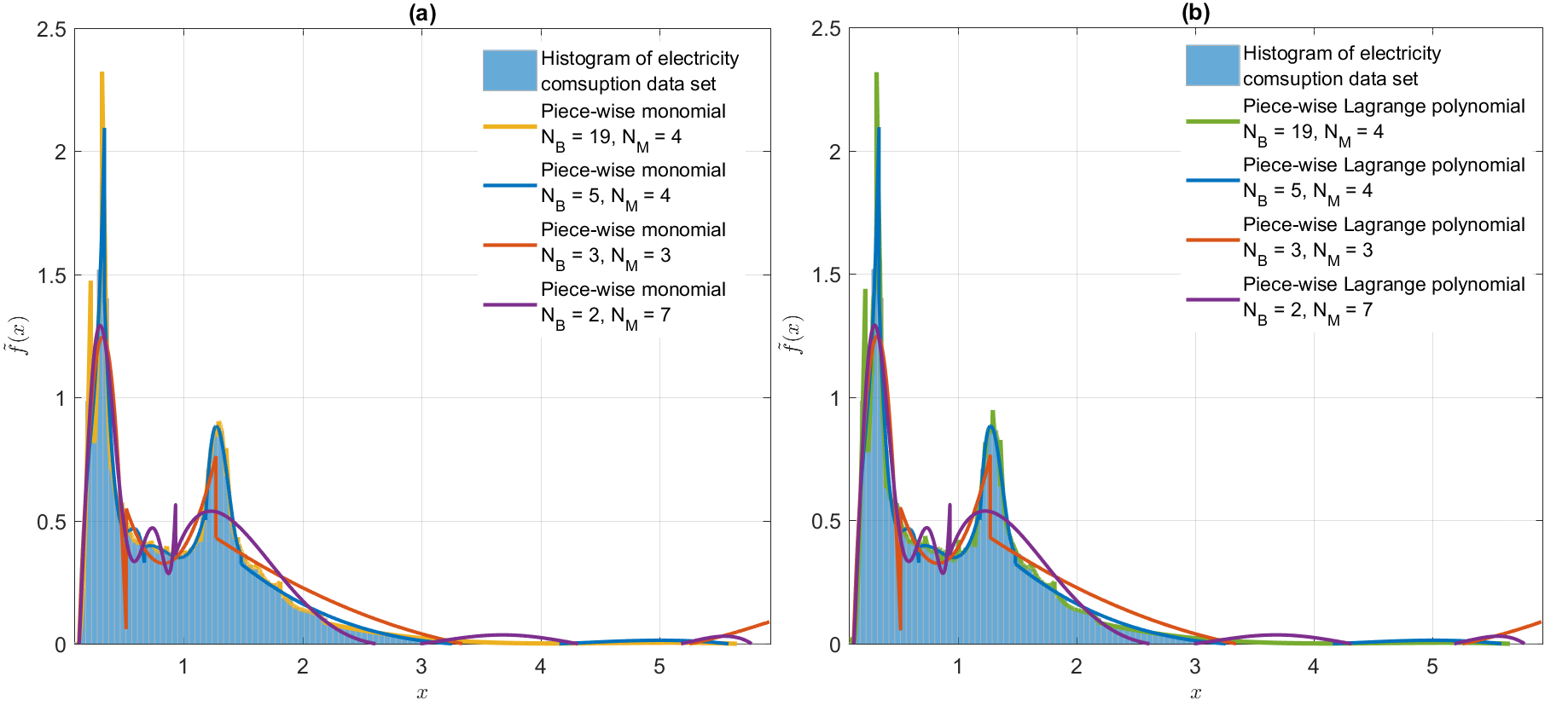}
    \caption{Sensitivity of the K--S test statistic to the number of bins $N_B$ and polynomial moment order $N_M$ for the household electricity consumption dataset, for both piecewise monomial and piecewise Lagrange polynomial approaches.}
    \label{fig:Electricity_sensitivity}
\end{figure}

\begin{figure}[H]
    \centering
    \includegraphics[width=1\linewidth]{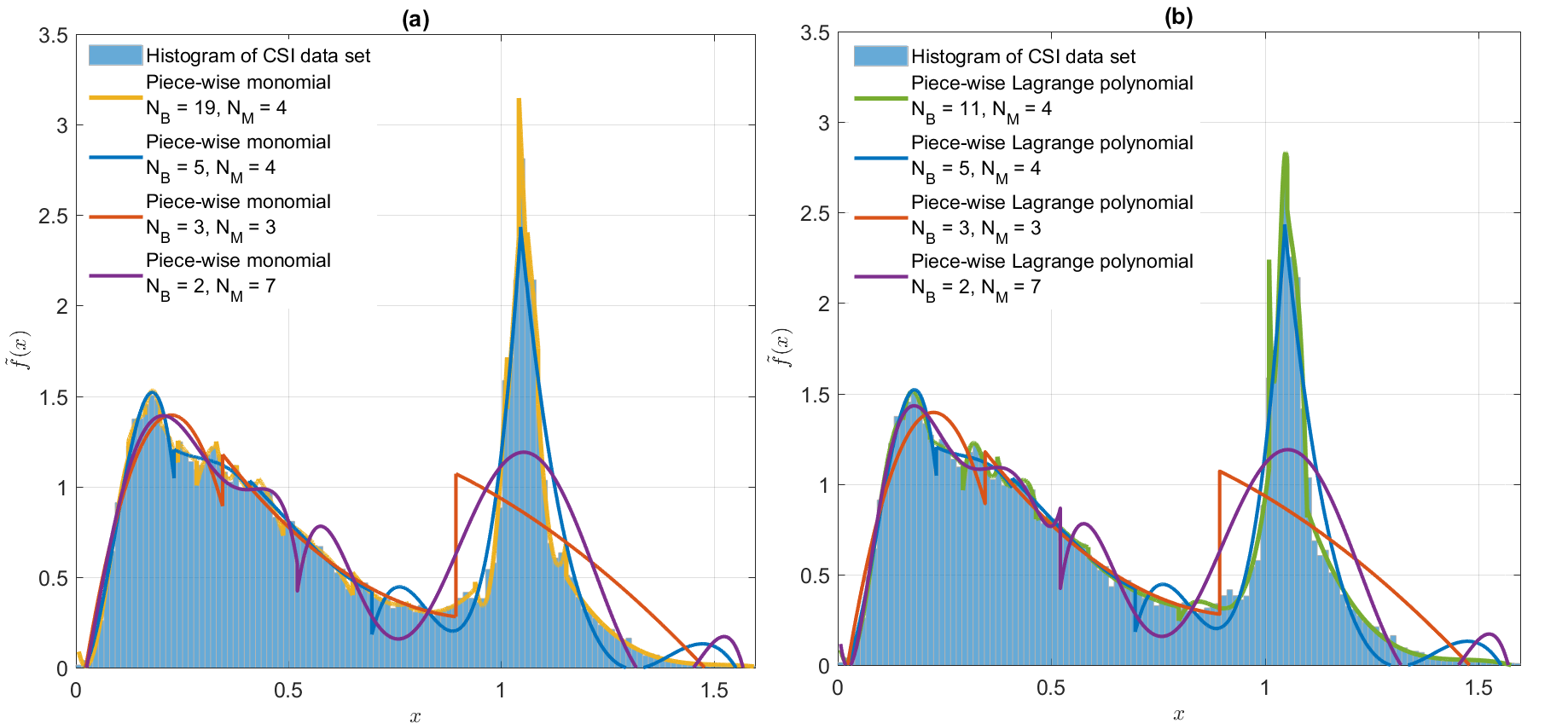}
    \caption{Sensitivity of the K--S test statistic to the number of bins $N_B$ and polynomial moment order $N_M$ for the solar irradiance dataset, for both piecewise monomial and piecewise Lagrange polynomial approaches.}
    \label{fig:CSI_sensitivity}
\end{figure}

Figures~\ref{fig:Electricity_sensitivity} and \ref{fig:CSI_sensitivity} display the sensitivity maps for the household electricity consumption and solar irradiance datasets, respectively, illustrating how the K--S test statistic varies across the $(N_B,N_M)$ grid. Tables~\ref{tab:KS_both_datasets} and the supplementary tables provide detailed numerical values.

Table \ref{tab:KS_both_datasets} presents a sensitivity analysis of the piecewise monomial and piecewise Lagrange polynomial approaches for the household electricity consumption and solar irradiance datasets. The first row in each section corresponds to the minimum K--S test statistic achieved, representing the optimal combination of bin number ($N_B$) and polynomial degree ($N_M$). For household electricity consumption, piecewise monomial approach achieves its lowest K--S test statistic value of 1.75e-03 at $N_B=19$ and $N_M=4$, while piecewise Lagrange approach reaches 1.96e-03 for the same configuration. This indicates that for relatively symmetric and smooth distributions, piecewise monomial approach provides a slightly better approximation. In contrast, for the more variable solar irradiance dataset, piecewise Lagrange polynomial approach attains the lowest K--S test statistic value of 1.92e-03 at $N_B=11$ and $N_M=4$, slightly outperforming piecewise monomial approach (2.13e-03), highlighting its advantage in capturing skewness, and localized variations.

The subsequent rows illustrate the effect of reducing the number of bins or changing the polynomial degree. As $N_B$ decreases or $N_M$ is modified away from the optimal values, the K--S test statistics increase, confirming that both the number of bins and the polynomial degree are important for achieving accurate PDF reconstruction. Overall, the sensitivity analysis demonstrates that while both piecewise approaches are robust, careful selection of $N_B$ and $N_M$ is crucial, with piecewise monomial favoring smoother, symmetric distributions and piecewise Lagrange polynomial better suited for more skewed or heterogeneous datasets.

\begin{table}[H]
\caption{Selected K--S test statistics for piecewise monomial and piecewise Lagrange polynomial reconstructions for the household electricity consumption and solar irradiance datasets.}
\label{tab:KS_both_datasets}
\centering
\small
\setlength{\tabcolsep}{4pt}
\begin{tabular}{ccc|ccc || ccc|ccc}
\toprule
\multicolumn{6}{c||}{\textbf{Household electricity consumption}} &
\multicolumn{6}{c}{\textbf{solar irradiance}} \\
\cmidrule(r){1-6}\cmidrule(l){7-12}
\multicolumn{3}{c}{\textbf{piecewise monomial}} &
\multicolumn{3}{c||}{\shortstack{\textbf{piecewise Lagrange}\\\textbf{polynomial}}} &
\multicolumn{3}{c}{\textbf{piecewise monomial}} &
\multicolumn{3}{c}{\shortstack{\textbf{piecewise Lagrange }\\\textbf{polynomial}}} \\
\cmidrule(r){1-3}\cmidrule(lr){4-6}\cmidrule(lr){7-9}\cmidrule(l){10-12}
$N_B$ & $N_M$ & $K\!-\!S$ &
$N_B$ & $N_M$ & $K\!-\!S$ &
$N_B$ & $N_M$ & $K\!-\!S$ &
$N_B$ & $N_M$ & $K\!-\!S$ \\
\midrule
19 & 4 & 1.75e-03 & 19 & 4 & 1.96e-03 &
19 & 4 & 2.13e-03 & 11 & 4 & 1.92e-03 \\
5  & 4 & 7.20e-03 & 5  & 4 & 7.20e-03 &
5  & 4 & 1.24e-02  & 5  & 4 & 1.24e-02  \\
3  & 3 & 4.97e-02  & 3  & 3 & 4.97e-02   &
3  & 3 & 5.69e-02  & 3  & 3 & 5.69e-02  \\
2  & 7 & 3.08e-02  & 2  & 7 & 3.08e-02  &
2  & 7 & 4.09e-02  & 2  & 7 & 4.09e-02  \\
\bottomrule
\end{tabular}

\end{table}

\section{Conclusion}

This work has introduced a piecewise polynomial with quantile-based binning framework for probability density reconstruction using the method of moments. By dividing the data domain into bins defined by empirical quantiles and fitting low-order polynomials locally in each bin, the approach mitigates the oscillatory behavior, tail instability, and sensitivity to support selection that often affect moment-based polynomial approximations.

Two formulations were investigated: piecewise monomial and piecewise Lagrange polynomials with Chebyshev nodes. In both cases, local coefficients are obtained by matching empirical moments within each bin, while the overall model complexity is controlled through the choice of bin count and polynomial degree. A grid search guided by the Kolmogorov--Smirnov statistic, together with a non-negativity feasibility check, provides a practical strategy for selecting these structural parameters.

Applications to benchmark distributions and real-world data sets in energy engineering show that the proposed piecewise approaches consistently improve the agreement with empirical distributions compared with monomial, Lagrange polynomial, and Spline Density approaches, and can perform competitively with Kernel Density Estimation. The findings indicate that quantile-based localization offers a flexible and robust pathway for moment-driven density reconstruction, particularly in settings where accurate tail behavior, heterogeneous and multimodal structures are of interest.

\section{Disclosure statement}\label{disclosure-statement}

The authors declare that they have no conflicts of interest.

\section{Data Availability Statement}\label{data-availability-statement}

Data supporting the findings of this study are available from the corresponding author upon reasonable request.

\phantomsection\label{supplementary-material}
\bigskip

\begin{center}

{\large\bf SUPPLEMENTARY MATERIAL}

\end{center}

\textbf{Title: Kolmogorov-Smirnov statistics and Goodness-of-fit measurement of piecewise monomial and piecewise Lagrange polynomial approaches with $N_B=1,...,19$ and $N_M=3,...,11$ for the household electricity consumption dataset}

\begin{table}[H]
\caption{Kolmogorov--Smirnov statistics $K-S$ of household electricity consumption data set as a function of the number of bins $N_B$ and the polynomial moment order $N_M$ for piecewise monomial approach}
    \label{tab:KS_NB_NM_matrix}
    \centering
    \scriptsize
    \begin{tabular}{c|ccccccccc}
        $N_B$ & $N_M=3$ & $N_M=4$ & $N_M=5$ & $N_M=6$ & $N_M=7$ & $N_M=8$ & $N_M=9$ & $N_M=10$ & $N_M=11$ \\
        \hline
        1  & 4.95e-02 & 5.55e-02 & 5.67e-02 & 6.10e-02 & 5.66e-02 & 5.62e-02 & 5.29e-02 & 4.80e-02 & 5.43e-02 \\
        2  & 5.38e-02 & 4.36e-02 & 4.01e-02 & 3.43e-02 & 3.08e-02 & 2.37e-02 & 1.88e-02 & 2.49e-02 & 1.88e-02 \\
        3  & 4.97-02 & 2.89-02 & 1.87e-02 & 1.29e-02 & 9.41e-03 & 9.68e-03 & 9.72e-03 & 9.75e-03 & 9.77e-03 \\
        4  & 2.85e-02 & 1.07e-02 & 8.45e-03 & 8.28e-03 & 8.18e-03 & 8.13e-03 &8.11e-03 & 8.11e-03 &8.11e-03 \\
        5  & 2.21e-02 & 7.20e-03 & 6.27e-03 & 6.21e-03 & 6.19e-03 & 6.19e-03 & 6.19e-03 & 6.19e-03 & 6.19e-03 \\
        6  & 1.84e-02 & 6.67e-03 & 6.25e-03 & 6.18e-03 & 6.16e-03 & 6.15e-03 & 6.16e-03 & 6.15e-03 & 6.15e-03 \\
        7  & 1.43e-02 & 6.41e-03 & 6.30e-03 & 6.26e-03 & 6.25e-03 & 6.24e-03 & 6.24e-03 & 6.24e-03 & 6.24e-03 \\
        8  & 1.27e-02 & 6.16e-03 & 6.11e-03 & 6.09e-03 & 6.09e-03 & 6.09e-03 & 6.09e-03 & 6.09e-03 & 6.09e-03 \\
        9  & 1.11e-02 & 5.41e-03 & 5.42e-03 & 5.43e-03 & 5.43e-03 & 5.43e-03 & 5.43e-03 & 5.43e-03 & 5.43e-03 \\
        10 & 9.64e-03 & 4.39e-03 & 4.43e-03 & 4.45e-03 & 4.45e-03 & 4.45e-03 & 4.45e-03 & 4.45e-03 & 4.45e-03 \\
        11 & 8.64e-03 & 3.36e-03 & 3.40e-03 & 3.42e-03 & 3.42e-03 & 3.42e-03 & 3.42e-03 & 3.42e-03 & 3.42e-03 \\
        12 & 8.06e-03 & 2.69e-03 & 2.48e-03 & 2.47e-03 & 2.46e-03 & 2.46e-03 & 2.46e-03 & 2.46e-03 & 2.46e-03 \\
        13 & 7.51e-03 & 3.12e-03 & 3.07e-03 & 3.05e-03 & 3.05e-03 & 3.05e-03 & 3.05e-03 & 3.05e-03 & 3.05e-03 \\
        14 & 7.07e-03 & 3.19e-03 & 3.17e-03 & 3.16e-03 & 3.16e-03 & 3.16e-03 & 3.16e-03 & 3.16e-03 & 3.16e-03 \\
        15 & 6.63e-03 & 2.80e-03 & 2.80e-03 & 2.80e-03 & 2.80e-03 & 2.80e-03 & 2.80e-03 & 2.80e-03 & 2.80e-03 \\
        16 & 6.25e-03 & 2.80e-03 & 2.80e-03 & 2.81e-03 & 2.81e-03 & 2.81e-03 & 2.81e-03 & 2.81e-03 & 2.81e-03 \\
        17 & 5.76e-03 & 2.26e-03 & 2.26e-03 & 2.26e-03 & 2.27e-03 & 2.27e-03 & 2.27e-03 & 2.27e-03 & 2.27e-03 \\
        18 & 5.42e-03 & 2.09e-03 & 2.10e-03 & 2.10e-03 & 2.10e-03 & 2.10e-03 & 2.10e-03 & 2.10e-03 & 2.10e-03 \\
        19 & 5.17e-03 & 1.75e-03 & 1.96e-03 & 1.96e-03 & 1.96e-03 & 1.96e-03 &1.96e-03 & 1.96e-03 & 1.96e-03 \\
    \end{tabular}
    
\end{table}

\begin{table}[H]
\caption{$GoF$ of household electricity consumption data set as a function of the number of bins $N_B$ and the polynomial moment order $N_M$ for piecewise monomial approach}
\centering
\scriptsize
\begin{tabular}{c|ccccccccc}
\hline
$N_B$ & $N_M=3$ & $N_M=4$ & $N_M=5$ & $N_M=6$ & $N_M=7$ & $N_M=8$ & $N_M=9$ & $N_M=10$ & $N_M=11$ \\
\hline
1&0.96 & 0.97 & 0.97 & 0.97 & 0.97 & 0.97 & 0.97 & 0.97 & 0.97 \\
2&0.95 & 0.97 & 0.98 & 0.98 & 0.98 & 0.98 & 0.99 & 0.99 & 0.99 \\
3&0.96 & 0.98 & 0.99 & 0.99 & 0.99 & 0.99 & 0.99 & 0.99 & 0.99 \\
4&0.98 & 0.99 & 0.99 & 0.99 & 0.99 & 0.99 & 0.99 & 0.99 & 0.99 \\
5&0.98 & 0.99 & 0.99 & 0.99 & 0.99 & 0.99 & 0.99 & 0.99 & 0.99 \\
6&0.99 & 0.99 & 0.99 & 0.99 & 0.99 & 0.99 & 0.99 & 0.99 & 0.99 \\
7&0.99 & 0.99 & 0.99 & 0.99 & 0.99 & 0.99 & 0.99 & 0.99 & 0.99 \\
8&0.99 & 0.99 & 0.99 & 0.99 & 0.99 & 0.99 & 0.99 & 0.99 & 0.99 \\
9&0.99 & 0.99 & 0.99 & 0.99 & 0.99 & 0.99 & 0.99 & 0.99 & 0.99 \\
10&0.99 & 0.99 & 0.99 & 0.99 & 0.99 & 0.99 & 0.99 & 0.99 & 0.99 \\
11&0.99 & 0.99 & 0.99 & 0.99 & 0.99 & 0.99 & 0.99 & 0.99 & 0.99 \\
12&0.99 & 0.99 & 0.99 & 0.99 & 0.99 & 0.99 & 0.99 & 0.99 & 0.99 \\
13&0.99 & 0.99 & 0.99 & 0.99 & 0.99 & 0.99 & 0.99 & 0.99 & 0.99 \\
14&0.99 & 0.99 & 0.99 & 0.99 & 0.99 & 0.99 & 0.99 & 0.99 & 0.99 \\
15&0.99 & 0.99 & 0.99 & 0.99 & 0.99 & 0.99 & 0.99 & 0.99 & 0.99 \\
16&0.99 & 0.99 & 0.99 & 0.99 & 0.99 & 0.99 & 0.99 & 0.99 & 0.99 \\
17&0.99 & 0.99 & 0.99 & 0.99 & 0.99 & 0.99 & 0.99 & 0.99 & 0.99 \\
18&0.99 & 0.99 & 0.99 & 0.99 & 0.99 & 0.99 & 0.99 & 0.99 & 0.99 \\
19&0.99 & 0.99 & 0.99 & 0.99 & 0.99 & 0.99 & 0.99 & 0.99 & 0.99 \\
\hline
\end{tabular}

\end{table}

\begin{table}[H]
\caption{Kolmogorov--Smirnov statistics $K$--$S$ of household electricity consumption data set as a function of the number of bins $N_B$ and the polynomial moment order $N_M$ for piecewise Lagrange polynomial approach}
    \label{tab:KS_NB_NM}
    \centering
    \scriptsize
    \begin{tabular}{c|ccccccccc}
        \hline
        $N_B$ & $N_M=3$ & $N_M=4$ & $N_M=5$ & $N_M=6$ & $N_M=7$ & $N_M=8$ & $N_M=9$ & $N_M=10$ & $N_M=11$ \\
        \hline
        1  & 4.95e-02 & 5.55e-02 & 5.67e-02 & 6.10e-02 & 5.66e-02 & 5.62e-02 & 5.29e-02 & 4.81e-02 & 4.30e-02 \\
        2  & 5.38e-02 & 4.36e-02 & 4.01e-02 & 3.43e-02 & 3.08e-02 & 2.37e-02 & 1.82e-02 & 1.73e-02 & 1.76e-02 \\
        3  & 4.97e-02 & 2.89e-02 & 1.87e-02 & 1.29e-02 & 9.41e-03 & 7.36e-03 & 5.97e-03 & 5.76e-03 & 5.47e-03 \\
        4  & 2.85e-02 & 1.07e-02 & 8.45e-03 & 5.47e-03 & 4.37e-03 & 3.26e-03 & 7.41e-03 & 1.40e-02 & 1.03e-01 \\
        5  & 2.21e-02 & 7.20e-03 & 5.48e-03 & 3.58e-03 & 3.47e-03 & 4.21e-02 & 8.75e-02 & 1.05e-01 & 4.77e-02 \\
        6  & 1.84e-02 & 6.67e-03 & 3.94e-03 & 3.50e-03 & 2.90e-03 & 2.81e-02 & 4.31e-02 & 1.08e-01 & 2.30e-01 \\
        7  & 1.43e-02 & 6.41e-03 & 3.33e-03 & 3.30e-03 & 3.18e-03 & 1.73e-02 & 8.28e-02 & 6.35e-02 & 6.25e-01 \\
        8  & 1.27e-02 & 6.16e-03 & 3.19e-03 & 1.47e-02 & 3.90e-02 & 1.71e-02 & 2.48e+00 & 9.07e+00 & 4.58e-01 \\
        9  & 1.11e-02 & 5.41e-03 & 3.45e-03 & 2.72e-03 & 3.48e-02 & 2.55e-02 & 3.77e-02 & 4.12e-01 & 4.33e-01 \\
        10 & 9.64e-03 & 4.39e-03 & 3.39e-03 & 1.91e-02 & 6.42e-02 & 3.88e-02 & 4.90e-01 & 8.74e-01 & 2.19e-01 \\
        11 & 8.65e-03 & 3.36e-03 & 3.09e-03 & 1.22e-02 & 7.93e-02 & 8.87e-01 & 3.67e-01 & 2.91e-01 & 4.24e+00 \\
        12 & 8.06e-03 & 2.69e-03 & 2.73e-03 & 1.02e-02 & 4.63e-01 & 2.37e+00 & 4.71e-01 & 1.32e-01 & 5.40e-01 \\
        13 & 7.51e-03 & 3.12e-03 & 2.57e-03 & 7.85e-03 & 2.39e+00 & 5.38e-01 & 5.78e-02 & 1.36e+00 & 2.27e-01 \\
        14 & 7.07e-03 & 3.19e-03 & 2.59e-03 & 3.29e-02 & 6.81e-02 & 1.38e-01 & 1.11e+00 & 7.50e-01 & 5.11e-01 \\
        15 & 6.63e-03 & 2.80e-03 & 2.64e-03 & 4.43e-02 & 7.02e-02 & 2.40e-01 & 6.28e-01 & 1.96e-01 & 5.33e-01 \\
        16 & 6.25e-03 & 2.80e-03 & 4.62e-03 & 4.23e-02 & 1.25e-01 & 7.49e-02 & 2.13e-01 & 1.05e-01 & 5.74e-01 \\
        17 & 5.77e-03 & 2.26e-03 & 2.35e-03 & 5.60e-02 & 1.12e-01 & 3.27e-02 & 2.82e-01 & 1.44e+00 & 3.25e-01 \\
        18 & 5.42e-03 & 2.09e-03 & 2.52e-03 & 9.00e-03 & 3.06e+00 & 1.85e-01 & 1.23e-01 & 3.94e-01 & 6.05e-01 \\
        19 & 5.17e-03 & 1.96e-03 & 2.38e-03 & 2.36e-02 & 5.52e-02 & 2.21e+01 & 1.36e-01 & 2.64e-01 & 4.19e-01 \\
        \hline
    \end{tabular}
    
\end{table}

\begin{table}[H]
\caption{$GoF$ of household electricity consumption data set as a function of the number of bins $N_B$ and the polynomial moment order $N_M$ for piecewise Lagrange polynomial approach}
\centering
\scriptsize
\begin{tabular}{c|ccccccccc}
\hline
$N_B$ & $N_M=3$ & $N_M=4$ & $N_M=5$ & $N_M=6$ & $N_M=7$ & $N_M=8$ & $N_M=9$ & $N_M=10$ & $N_M=11$ \\
        \hline
1&0.96 & 0.97 & 0.97 & 0.97 & 0.97 & 0.97 & 0.97 & 0.97 & 0.98 \\
2&0.95 & 0.97 & 0.98 & 0.98 & 0.98 & 0.99 & 0.99 & 0.99 & 0.99 \\
3&0.96 & 0.98 & 0.99 & 0.99 & 0.99 & 0.99 & 0.99 & 0.99 & 0.99 \\
4&0.98 & 0.99 & 0.99 & 0.99 & 0.99 & 0.99 & 0.99 & 0.99 & 0.94 \\
5&0.98 & 0.99 & 0.99 & 0.99 & 0.99 & 0.97 & 0.95 & 0.95 & 0.97 \\
6&0.99 & 0.99 & 0.99 & 0.99 & 0.99 & 0.98 & 0.98 & 0.94 & 0.88 \\
7&0.99 & 0.99 & 0.99 & 0.99 & 0.99 & 0.99 & 0.96 & 0.96 & 0.74 \\
8&0.99 & 0.99 & 0.99 & 0.99 & 0.98 & 0.99 & 0.08 & -0.41 & 0.83 \\
9&0.99 & 0.99 & 0.99 & 0.99 & 0.98 & 0.98 & 0.97 & 0.84 & 0.49 \\
10&0.99 & 0.99 & 0.99 & 0.99 & 0.97 & 0.98 & 0.76 & 0.58 & 0.81 \\
11&0.99 & 0.99 & 0.99 & 0.99 & 0.97 & 0.71 & 0.87 & 0.86 & 4.07 \\
12&0.99 & 0.99 & 0.99 & 0.99 & 0.82 & 0.19 & 0.85 & 0.90 & 0.62 \\
13&0.99 & 0.99 & 0.99 & 0.99 & 0.30 & 0.81 & 0.97 & -1.97 & 0.79 \\
14&0.99 & 0.99 & 0.99 & 0.98 & 0.97 & 0.95 & 0.64 & 0.59 & 0.27 \\
15&0.99 & 0.99 & 0.99 & 0.98 & 0.97 & 0.92 & 0.79 & 0.81 & 0.09 \\
16&0.99 & 0.99 & 0.99 & 0.98 & 0.94 & 0.97 & 0.91 & 0.91 & -0.07 \\
17&0.99 & 0.99 & 0.99 & 0.98 & 0.96 & 0.98 & 0.88 & 0.28 & 0.51 \\
18&0.99 & 0.99 & 0.99 & 0.99 & 0.14 & 0.94 & 0.89 & 0.40 & 0.73 \\
19&0.99 & 0.99 & 0.99 & 0.99 & 0.97 & -5.05 & 0.88 & 0.67 & 0.48 \\
\hline
\end{tabular}
\end{table}
\newpage

\textbf{Title: Kolmogorov-Smirnov statistics and Goodness-of-fit measurement of piecewise monomial and piecewise Lagrange polynomial approaches with $N_B=1,...,19$ and $N_M=3,...,11$ for the solar irradiance dataset}

\begin{table}[H]
\caption{Kolmogorov--Smirnov statistics $K-S$ of solar irradiance data set as a function of the number of bins $N_B$ and the polynomial moment order $N_M$ for piecewise monomial approach}
\centering
\scriptsize
\begin{tabular}{c|ccccccccc}
        $N_B$ & $N_M=3$ & $N_M=4$ & $N_M=5$ & $N_M=6$ & $N_M=7$ & $N_M=8$ & $N_M=9$ & $N_M=10$ & $N_M=11$ \\
        \hline
1&1.03e-01 & 1.00e-01 & 7.32e-02 & 6.13e-02 & 5.42e-02 & 5.68e-02 & 4.48e-02 & 4.31e-02 & 4.50e-02 \\
2&7.47e-02 & 6.92e-02 & 5.24e-02 & 5.20e-02 & 4.09e-02 & 3.96e-02 & 4.13e-02 & 4.16e-02 & 4.07e-02 \\
3&5.69e-02 & 4.72e-02 & 2.92e-02 & 3.00e-02 & 2.96e-02 & 2.98e-02 & 2.97e-02 & 2.83e-02 & 2.82e-02 \\
4&2.78e-02 & 1.33e-02 & 1.08e-02 & 1.08e-02 & 1.08e-02 & 1.08e-02 & 1.08e-02 & 1.07e-02 & 1.06e-02 \\
5&2.63e-02 & 1.24e-02 & 6.74e-03 & 7.07e-03 & 7.39e-03 & 4.40e-03 & 4.75e-03 & 5.12e-03 & 5.50e-03 \\
6&2.06e-02 & 1.12e-02 & 6.09e-03 & 6.55e-03 & 7.07e-03 & 4.22e-03 & 4.46e-03 & 4.71e-03 & 4.96e-03 \\
7&1.59e-02 & 9.63e-03 & 5.90e-03 & 6.29e-03 & 6.70e-03 & 4.61e-03 & 4.77e-03 & 4.91e-03 & 5.03e-03 \\
8&1.18e-02 & 7.12e-03 & 4.80e-03 & 5.04e-03 & 5.31e-03 & 4.01e-03 & 4.07e-03 & 4.12e-03 & 4.14e-03 \\
9&7.80e-03 & 4.67e-03 & 2.90e-03 & 3.08e-03 & 3.27e-03 & 3.48e-03 & 2.66e-03 & 2.67e-03 & 2.67e-03 \\
10&4.66e-03 & 2.87e-03 & 2.82e-03 & 2.81e-03 & 2.81e-03 & 2.81e-03 & 2.71e-03 & 2.70e-03 & 2.69e-03 \\
11&3.17e-03 & 2.44e-03 & 2.44e-03 & 2.44e-03 & 2.43e-03 & 2.43e-03 & 2.43e-03 & 2.42e-03 & 2.42e-03 \\
12&2.78e-03 & 2.70e-03 & 2.63e-03 & 2.96e-03 & 2.96e-03 & 2.96e-03 & 2.95e-03 & 2.95e-03 & 2.94e-03 \\
13&3.29e-03 & 2.80e-03 & 2.81e-03 & 2.81e-03 & 2.81e-03 & 2.81e-03 & 2.68e-03 & 2.68e-03 & 2.68e-03 \\
14&2.56e-03 & 2.56e-03 & 2.56e-03 & 2.56e-03 & 2.57e-03 & 2.57e-03 & 2.57e-03 & 2.57e-03 & 2.58e-03 \\
15&2.39e-03 & 2.35e-03 & 2.36e-03 & 2.36e-03 & 2.36e-03 & 2.36e-03 & 2.36e-03 & 2.36e-03 & 2.36e-03 \\
16&2.81e-03 & 2.81e-03 & 2.81e-03 & 2.82e-03 & 2.82e-03 & 2.82e-03 & 2.33e-03 & 2.31e-03 & 2.30e-03 \\
17&2.63e-03 & 2.63e-03 & 2.63e-03 & 2.63e-03 & 2.63e-03 & 2.63e-03 & 2.63e-03 & 2.64e-03 & 2.64e-03 \\
18&2.29e-03 & 2.34e-03 & 2.34e-03 & 2.34e-03 & 2.34e-03 & 2.34e-03 & 2.34e-03 & 2.34e-03 & 2.34e-03 \\
19&2.15e-03 & 2.13e-03 & 2.16e-03 & 2.16e-03 & 2.16e-03 & 2.16e-03 & 2.16e-03 & 2.16e-03 & 2.16e-03 \\
\end{tabular}
\end{table}

\begin{table}[H]
\caption{$GoF$ of solar irradiance data set as a function of the number of bins $N_B$ and the polynomial moment order $N_M$ for piecewise monomial approach}
\centering
\scriptsize
\begin{tabular}{c|ccccccccc}
$N_B$ & $N_M=3$ & $N_M=4$ & $N_M=5$ & $N_M=6$ & $N_M=7$ & $N_M=8$ & $N_M=9$ & $N_M=10$ & $N_M=11$ \\
\hline
1  & 0.92 & 0.93 & 0.94 & 0.96 & 0.96 & 0.96 & 0.97 & 0.97 & 0.97 \\
2  & 0.95 & 0.95 & 0.97 & 0.97 & 0.97 & 0.97 & 0.97 & 0.98 & 0.98 \\
3  & 0.96 & 0.97 & 0.98 & 0.98 & 0.98 & 0.98 & 0.98 & 0.98 & 0.98 \\
4  & 0.98 & 0.99 & 0.99 & 0.99 & 0.99 & 0.99 & 0.99 & 0.99 & 0.99 \\
5  & 0.98 & 0.99 & 0.99 & 0.99 & 0.99 & 0.99 & 0.99 & 0.99 & 0.99 \\
6  & 0.99 & 0.99 & 0.99 & 0.99 & 0.99 & 0.99 & 0.99 & 0.99 & 0.99 \\
7  & 0.99 & 0.99 & 0.99 & 0.99 & 0.99 & 0.99 & 0.99 & 0.99 & 0.99 \\
8  & 0.99 & 0.99 & 0.99 & 0.99 & 0.99 & 0.99 & 0.99 & 0.99 & 0.99 \\
9  & 0.99 & 0.99 & 0.99 & 0.99 & 0.99 & 0.99 & 0.99 & 0.99 & 0.99 \\
10 & 0.99 & 0.99 & 0.99 & 0.99 & 0.99 & 0.99 & 0.99 & 0.99 & 0.99 \\
11 & 0.99 & 0.99 & 0.99 & 0.99 & 0.99 & 0.99 & 0.99 & 0.99 & 0.99 \\
12 & 0.99 & 0.99 & 0.99 & 0.99 & 0.99 & 0.99 & 0.99 & 0.99 & 0.99 \\
13 & 0.99 & 0.99 & 0.99 & 0.99 & 0.99 & 0.99 & 0.99 & 0.99 & 0.99 \\
14 & 0.99 & 0.99 & 0.99 & 0.99 & 0.99 & 0.99 & 0.99 & 0.99 & 0.99 \\
15 & 0.99 & 0.99 & 0.99 & 0.99 & 0.99 & 0.99 & 0.99 & 0.99 & 0.99 \\
16 & 0.99 & 0.99 & 0.99 & 0.99 & 0.99 & 0.99 & 0.99 & 0.99 & 0.99 \\
17 & 0.99 & 0.99 & 0.99 & 0.99 & 0.99 & 0.99 & 0.99 & 0.99 & 0.99 \\
18 & 0.99 & 0.99 & 0.99 & 0.99 & 0.99 & 0.99 & 0.99 & 0.99 & 0.99 \\
19 & 0.99 & 0.99 & 0.99 & 0.99 & 0.99 & 0.99 & 0.99 & 0.99 & 0.99 \\
\end{tabular}
\end{table}

\begin{table}[H]
\caption{Kolmogorov--Smirnov statistics $K-S$ of solar irradiance data set as a function of the number of bins $N_B$ and the polynomial moment order $N_M$ for piecewise Lagrange polynomial approach}
\label{tab:NB_NM_scientific}
\centering
\scriptsize
\begin{tabular}{c|ccccccccc}
$N_B$ & $N_M=3$ & $N_M=4$ & $N_M=5$ & $N_M=6$ & $N_M=7$ & $N_M=8$ & $N_M=9$ & $N_M=10$ & $N_M=11$ \\
\hline
1  & 1.03e-01 & 1.00e-01 & 7.32e-02 & 6.13e-02 & 5.42e-02 & 5.68e-02 & 4.48e-02 & 4.36e-02 & 3.84e-02 \\
2  & 7.47e-02 & 6.92e-02 & 5.24e-02 & 5.20e-02 & 4.09e-02 & 4.09e-02 & 3.23e-02 & 3.28e-02 & 3.10e-02 \\
3  & 5.69e-02 & 4.72e-02 & 2.92e-02 & 2.91e-02 & 2.68e-02 & 1.78e-02 & 1.46e-02 & 2.22e-02 & 3.32e-02 \\
4  & 2.78e-02 & 1.33e-02 & 1.08e-02 & 1.09e-02 & 8.82e-03 & 6.42e-03 & 2.50e-02 & 3.39e-02 & 7.43e-02 \\
5  & 2.63e-02 & 1.24e-02 & 6.75e-03 & 3.80e-03 & 3.30e-03 & 7.13e-03 & 3.65e-02 & 1.34e+00 & 1.07e+00 \\
6  & 2.06e-02 & 1.12e-02 & 6.09e-03 & 5.73e-03 & 1.95e-02 & 1.00e-02 & 2.37e-01 & 1.12e-01 & 1.91e+00 \\
7  & 1.59e-02 & 9.63e-03 & 5.90e-03 & 1.75e-02 & 1.44e-01 & 1.89e-02 & 7.69e-02 & 2.33e-01 & 4.15e-01 \\
8  & 1.18e-02 & 7.12e-03 & 4.89e-03 & 7.41e-03 & 7.16e-02 & 3.03e-02 & 9.25e-02 & 1.51e-01 & 1.19e+00 \\
9  & 7.80e-03 & 4.67e-03 & 3.05e-03 & 4.22e-02 & 2.63e-01 & 7.81e-01 & 3.69e+00 & 1.25e-01 & 4.59e-01 \\
10 & 4.66e-03 & 2.87e-03 & 4.84e-03 & 1.42e-02 & 6.54e-02 & 1.06e-01 & 1.86e-01 & 1.03e+00 & 5.53e-01 \\
11 & 3.17e-03 & 1.92e-03 & 2.36e-03 & 3.98e-02 & 5.40e-02 & 4.12e-02 & 4.23e-01 & 2.10e-01 & 5.41e-01 \\
12 & 2.78e-03 & 2.97e-03 & 2.46e-03 & 8.60e-02 & 2.43e-02 & 5.06e-02 & 2.93e-01 & 4.96e-01 & 3.56e-01 \\
13 & 2.80e-03 & 2.68e-03 & 2.59e-03 & 1.16e-02 & 3.31e-02 & 6.44e-02 & 1.35e-01 & 1.63e-01 & 2.56e+00 \\
14 & 2.56e-03 & 2.35e-03 & 2.23e-03 & 2.06e-01 & 8.65e-02 & 9.01e-02 & 1.95e-01 & 1.86e-01 & 1.97e+00 \\
15 & 2.39e-03 & 2.19e-03 & 3.49e-02 & 1.03e-02 & 7.62e-02 & 2.91e-01 & 2.00e-01 & 4.69e-01 & 9.18e-01 \\
16 & 2.81e-03 & 2.40e-03 & 5.52e-03 & 4.39e+00 & 7.51e-02 & 1.70e+00 & 5.21e-01 & 5.53e-01 & 1.33e+00 \\
17 & 2.63e-03 & 2.61e-03 & 1.49e-02 & 1.93e-02 & 6.73e-02 & 3.48e-01 & 3.15e+00 & 6.94e-01 & 2.53e-01 \\
18 & 2.34e-03 & 2.37e-03 & 5.30e-03 & 1.86e-02 & 4.98e-02 & 1.51e-01 & 7.60e-02 & 4.20e-01 & 5.93e-01 \\
19 & 2.16e-03 & 1.98e-03 & 2.48e-02 & 9.38e-01 & 3.64e-02 & 4.81e-01 & 1.52e-01 & 3.50e+00 & 2.56e-01 \\
\end{tabular}
\end{table}

\begin{table}[H]
\caption{$GoF$ of solar irradiance data set as a function of the number of bins $N_B$ and the polynomial moment order $N_M$ for piecewise Lagrange polynomial approach}
\label{tab:gof_trunc}
\centering
\scriptsize
\begin{tabular}{c|ccccccccc}
$N_B$ & $N_M=3$ & $N_M=4$ & $N_M=5$ & $N_M=6$ & $N_M=7$ & $N_M=8$ & $N_M=9$ & $N_M=10$ & $N_M=11$ \\
\hline
1 & 0.92 & 0.93 & 0.94 & 0.96 & 0.96 & 0.96 & 0.97 & 0.97 & 0.98 \\
2 & 0.95 & 0.95 & 0.97 & 0.97 & 0.98 & 0.98 & 0.98 & 0.98 & 0.98 \\
3 & 0.96 & 0.97 & 0.98 & 0.98 & 0.98 & 0.99 & 0.99 & 0.98 & 0.98 \\
4 & 0.98 & 0.99 & 0.99 & 0.99 & 0.99 & 0.99 & 0.98 & 0.98 & 0.95 \\
5 & 0.98 & 0.99 & 0.99 & 0.99 & 0.99 & 0.99 & 0.97 & 0.35 & 0.50 \\
6 & 0.99 & 0.99 & 0.99 & 0.99 & 0.99 & 0.99 & 0.90 & 0.95 & 0.21 \\
7 & 0.99 & 0.99 & 0.99 & 0.99 & 0.94 & 0.99 & 0.95 & 0.86 & 0.64 \\
8 & 0.99 & 0.99 & 0.99 & 0.99 & 0.97 & 0.98 & 0.92 & 0.91 & 0.50 \\
9 & 0.99 & 0.99 & 0.99 & 0.98 & 0.91 & 0.73 & -0.26 & 0.90 & 0.74 \\
10 & 0.99 & 0.99 & 0.99 & 0.99 & 0.97 & 0.95 & 0.87 & 0.34 & 0.79 \\
11 & 0.99 & 0.99 & 0.99 & 0.98 & 0.98 & 0.98 & 0.62 & 0.83 & 0.69 \\
12 & 0.99 & 0.99 & 0.99 & 0.96 & 0.98 & 0.97 & 0.89 & 0.84 & 0.68 \\
13 & 0.99 & 0.99 & 0.99 & 0.99 & 0.98 & 0.97 & 0.90 & 0.90 & 0.43 \\
14 & 0.99 & 0.99 & 0.99 & 0.93 & 0.97 & 0.96 & 0.89 & 0.89 & 0.45 \\
15 & 0.99 & 0.99 & 0.98 & 0.99 & 0.97 & 0.91 & 0.83 & 0.75 & -0.49 \\
16 & 0.99 & 0.99 & 0.99 & -0.30 & 0.97 & -0.32 & 0.80 & 0.29 & -11.41 \\
17 & 0.99 & 0.99 & 0.99 & 0.99 & 0.98 & 0.90 & 0.01 & 0.66 & 0.83 \\
18 & 0.99 & 0.99 & 0.99 & 0.99 & 0.98 & 0.91 & 0.94 & 0.81 & 0.10 \\
19 & 0.99 & 0.99 & 0.99 & 0.72 & 0.98 & 0.72 & 0.92 & 0.08 & 0.79 \\
\end{tabular}
\end{table}

\bibliography{References}

@article{Yu2023,
  author    = {Yue Yu and Pavel Loskot},
  title     = {Polynomial Distributions and Transformations},
  journal   = {Mathematics},
  year      = {2023},
  volume    = {11},
  number    = {4},
  pages     = {985},
  doi       = {10.3390/math11040985},
  publisher = {MDPI}
}

@article{Chan2013,
  author  = {Siu On Chan and Ilias Diakonikolas and Rocco A. Servedio and Xiao},
  title   = {Efficient Density Estimation via Piecewise Polynomial Approximations},
  journal = {arXiv preprint arXiv:1305.3207},
  year    = {2013}
}

@article{papp2014shape,
  title={Shape-constrained estimation using nonnegative splines},
  author={Papp, D{\'a}vid and Alizadeh, Farid},
  journal={Journal of Computational and graphical Statistics},
  volume={23},
  number={1},
  pages={211--231},
  year={2014},
  publisher={Taylor \& Francis}
}

@article{Freedman1981,
  author    = {David Freedman and Persi Diaconis},
  title     = {On the Histogram as a Density Estimator: L2 Theory},
  journal   = {Zeitschrift f{\"u}r Wahrscheinlichkeitstheorie und Verwandte Gebiete},
  year      = {1981},
  volume    = {57},
  number    = {4},
  pages     = {453--476},
  doi       = {10.1007/BF01025868}
}

@article{Abdous2001,
  author    = {Abdous, Belkacem and Bensaid, El Mostafa},
  title     = {Multivariate Local Polynomial Fitting for a Probability Distribution Function and Its Partial Derivatives},
  journal   = {Journal of Nonparametric Statistics},
  year      = {2001},
  volume    = {13},
  number    = {1},
  pages     = {77--94},
  doi       = {10.1080/10485250108832867}
}

@article{munkhammar2017polynomial,
  title={Polynomial probability distribution estimation using the method of moments},
  author={Munkhammar, Joakim and Mattsson, Lars and Ryd{\'e}n, Jesper},
  journal={PloS one},
  volume={12},
  number={4},
  pages={e0174573},
  year={2017},
  publisher={Public Library of Science San Francisco, CA USA}
}

@article{turan2024polynomial,
  title={Polynomial approaches in improving accuracy of probability distribution estimation using the method of moments},
  author={Turan, Meltem and Munkhammar, Joakim and Dutta, Abhishek},
  journal={Journal of Chemical Technology \& Biotechnology},
  volume={99},
  number={5},
  pages={1056--1068},
  year={2024},
  publisher={Wiley Online Library}
}

@article{han2019kernel,
  title={Kernel density estimation model for wind speed probability distribution with applicability to wind energy assessment in China},
  author={Han, Qinkai and Ma, Sai and Wang, Tianyang and Chu, Fulei},
  journal={Renewable and Sustainable Energy Reviews},
  volume={115},
  pages={109387},
  year={2019},
  publisher={Elsevier}
}

@article{munkhammar2021very,
  title={Very short term load forecasting of residential electricity consumption using the Markov-chain mixture distribution (MCM) model},
  author={Munkhammar, Joakim and van der Meer, Dennis and Wid{\'e}n, Joakim},
  journal={Applied Energy},
  volume={282},
  pages={116180},
  year={2021},
  publisher={Elsevier}
}

@article{munkhammar2018markov,
  title={A Markov-chain probability distribution mixture approach to the clear-sky index},
  author={Munkhammar, Joakim and Wid{\'e}n, Joakim},
  journal={Solar Energy},
  volume={170},
  pages={174--183},
  year={2018},
  publisher={Elsevier}
}

@article{pavlides2022non,
  title={Non-parametric kernel-based estimation and simulation of precipitation amount},
  author={Pavlides, Andrew and Agou, Vasiliki D and Hristopulos, Dionissios T},
  journal={Journal of Hydrology},
  volume={612},
  pages={127988},
  year={2022},
  publisher={Elsevier}
}

@article{reyes2017bandwidth,
  title={Bandwidth selection in kernel density estimation for interval-grouped data},
  author={Reyes, Miguel and Francisco-Fern{\'a}ndez, Mario and Cao, Ricardo},
  journal={Test},
  volume={26},
  number={3},
  pages={527--545},
  year={2017},
  publisher={Springer}
}

@article{hsu2021moon,
  title={Moon image segmentation with a new mixture histogram model},
  author={Hsu, Chih-Yu and Shao, Lu-Jiao and Tseng, Kuo-Kun and Huang, Wan-Ting},
  journal={Enterprise Information Systems},
  volume={15},
  number={8},
  pages={1046--1069},
  year={2021},
  publisher={Taylor \& Francis}
}

@article{dhal2021histogram,
  title={Histogram Equalization Variants as Optimization Problems: A Review: KG Dhal et al.},
  author={Dhal, Krishna Gopal and Das, Arunita and Ray, Swarnajit and G{\'a}lvez, Jorge and Das, Sanjoy},
  journal={Archives of Computational Methods in Engineering},
  volume={28},
  number={3},
  pages={1471--1496},
  year={2021},
  publisher={Springer}
}

@article{kirkby2023spline,
  title={Spline local basis methods for nonparametric density estimation},
  author={Kirkby, J Lars and Leitao, {\'A}lvaro and Nguyen, Duy},
  journal={Statistic Surveys},
  volume={17},
  pages={75--118},
  year={2023},
  publisher={The American Statistical Association, the Bernoulli Society, the Institute~…}
}

@article{jain2022load,
  title={Load forecasting and risk assessment for energy market with renewable based distributed generation},
  author={Jain, Ritu and Mahajan, Vasundhara},
  journal={Renewable Energy Focus},
  volume={42},
  pages={190--205},
  year={2022},
  publisher={Elsevier}
}

@article{li2022integrated,
  title={An integrated missing-data tolerant model for probabilistic PV power generation forecasting},
  author={Li, Qiaoqiao and Xu, Yan and Chew, Benjamin Si Hao and Ding, Hongyuan and Zhao, Guopeng},
  journal={IEEE Transactions on Power Systems},
  volume={37},
  number={6},
  pages={4447--4459},
  year={2022},
  publisher={IEEE}
}

@article{gu2021forecasting,
  title={Forecasting and uncertainty analysis of day-ahead photovoltaic power using a novel forecasting method},
  author={Gu, Bo and Shen, Huiqiang and Lei, Xiaohui and Hu, Hao and Liu, Xinyu},
  journal={Applied Energy},
  volume={299},
  pages={117291},
  year={2021},
  publisher={Elsevier}
}

@article{shi2021wind,
  title={Wind speed distributions used in wind energy assessment: a review},
  author={Shi, Huanyu and Dong, Zhibao and Xiao, Nan and Huang, Qinni},
  journal={Frontiers in Energy Research},
  volume={9},
  pages={769920},
  year={2021},
  publisher={Frontiers Media SA}
}

@article{alibrandi2018kernel,
  title={Kernel density maximum entropy method with generalized moments for evaluating probability distributions, including tails, from a small sample of data},
  author={Alibrandi, Umberto and Mosalam, Khalid M},
  journal={International Journal for Numerical Methods in Engineering},
  volume={113},
  number={13},
  pages={1904--1928},
  year={2018},
  publisher={Wiley Online Library}
}

@article{valiant2017estimating,
  title={Estimating the unseen: improved estimators for entropy and other properties},
  author={Valiant, Gregory and Valiant, Paul},
  journal={Journal of the ACM (JACM)},
  volume={64},
  number={6},
  pages={1--41},
  year={2017},
  publisher={ACM New York, NY, USA}
}

@article{biggs2002modelling,
  title={Modelling activated sludge flocculation using population balances},
  author={Biggs, CA and Lant, PA},
  journal={Powder Technology},
  volume={124},
  number={3},
  pages={201--211},
  year={2002},
  publisher={Elsevier}
}

@article{gavaskar2018fast,
  title={Fast adaptive bilateral filtering},
  author={Gavaskar, Ruturaj G and Chaudhury, Kunal N},
  journal={IEEE transactions on Image Processing},
  volume={28},
  number={2},
  pages={779--790},
  year={2018},
  publisher={IEEE}
}

@book{conover1999practical,
  title={Practical nonparametric statistics},
  author={Conover, William Jay},
  year={1999},
  publisher={john wiley \& sons}
}

@article{munkhammar2014characterizing,
  title={Characterizing probability density distributions for household electricity load profiles from high-resolution electricity use data},
  author={Munkhammar, Joakim and Ryd{\'e}n, Jesper and Wid{\'e}n, Joakim},
  journal={Applied Energy},
  volume={135},
  pages={382--390},
  year={2014},
  publisher={Elsevier}
}

@book{silverman2018density,
  title={Density estimation for statistics and data analysis},
  author={Silverman, Bernard W},
  year={2018},
  publisher={Routledge}
}

@misc{SMHI2008,
title = {Swedish Meteorological and Hydrological Institute},
  author = {{Swedish Meteorological and Hydrological Institute}},
  year = {2008},
  note = {Norrköping, Sweden. Available at: http://www.smhi.se/en}
}

\end{document}